\documentclass{amsart}
 \usepackage[foot]{amsaddr}

\usepackage{amsthm}
\usepackage{amssymb}
\usepackage{amsmath}
\usepackage{amsfonts}
\usepackage{amsrefs}
\usepackage{textcomp}
\usepackage[T1]{fontenc}
\usepackage[utf8x]{inputenc}
\usepackage{textcomp}
\usepackage{wasysym}
\usepackage{stmaryrd}
\usepackage{esint}
\usepackage[all]{xy}
\usepackage{graphicx}
\usepackage{bbm}

\newtheorem{tw}{Theorem}[section]
\newtheorem{dfn}[tw]{Definition}
\newtheorem{uw}[tw]{Remark}
\newtheorem{prz}[tw]{Example}
\newtheorem{lem}[tw]{Lemma}
\newtheorem{stw}[tw]{Proposition}
\newtheorem{wn}[tw]{Corollary}
\newtheorem{hip}{Question}

\newtheorem*{ozn}{Notation}
\newtheorem*{dd}{Proof}
\newtheorem*{roz}{Rozwiązanie}
\newtheorem*{ak}{Acknowledgements}

\let\olddfn\dfn
\renewcommand{\dfn}{\olddfn\normalfont}

\let\oldozn\ozn
\renewcommand{\ozn}{\oldozn\normalfont}

\let\oldlem\lem
\renewcommand{\lem}{\oldlem\normalfont}

\let\oldstw\stw
\renewcommand{\stw}{\oldstw\normalfont}

\let\olduw\uw
\renewcommand{\uw}{\olduw\normalfont}

\let\oldwn\wn
\renewcommand{\wn}{\oldwn\normalfont}

\let\oldprz\prz
\renewcommand{\prz}{\oldprz\normalfont}

\let\olddd\dd
\renewcommand{\dd}{\olddd\normalfont}

\let\oldroz\roz
\renewcommand{\roz}{\oldroz\normalfont}

\let\oldak\ak
\renewcommand{\ak}{\oldak\normalfont}

\let\oldhip\hip
\renewcommand{\hip}{\oldhip\normalfont}

\def\Xint#1{\mathchoice
   {\XXint\displaystyle\textstyle{#1}}%
   {\XXint\textstyle\scriptstyle{#1}}%
   {\XXint\scriptstyle\scriptscriptstyle{#1}}%
   {\XXint\scriptscriptstyle\scriptscriptstyle{#1}}%
   \!\int}
\def\XXint#1#2#3{{\setbox0=\hbox{$#1{#2#3}{\int}$}
     \vcenter{\hbox{$#2#3$}}\kern-.5\wd0}}

\def\dashint{\Xint-}

\newcommand{\twopartdef}[4]
{
\left\{
		\begin{array}{ll}
			#1 & #2 \\
			#3 & #4
		\end{array}
	\right.
}

\newcommand{\threepartdef}[6]
{
	\left\{
		\begin{array}{lll}
			#1 & #2 \\
			#3 & #4 \\
			#5 & #6
		\end{array}
	\right.
}

\usepackage{geometry}
\author{Wojciech G\'{o}rny}

\address{W. G\'{o}rny: Faculty of Mathematics, Informatics and Mechanics, University of Warsaw, Warsaw, Poland.}

\email{w.gorny@mimuw.edu.pl}


\title{Existence of minimisers in the least gradient problem for general boundary data}


\begin{document}

\begin{abstract}
We study existence of minimisers to the least gradient problem on a strictly convex domain in two settings. On a bounded domain, we allow the boundary data to be discontinuous and prove existence of minimisers in terms of the Hausdorff measure of the discontinuity set. Later, we allow the domain to be unbounded, prove existence of minimisers and study their properties in terms of the regularity of boundary data and the shape of the domain. 
\end{abstract}

\maketitle

\section{Introduction}

The least gradient problem, studied extensively since the pioneering work of Sternberg-Williams-Ziemer, \cite{SWZ}, is the problem of minimalisation

\begin{equation}\label{problem}\tag{LGP}
\min \{ \int_\Omega |Du|, \quad u \in BV(\Omega), \quad u|_{\partial\Omega} = f  \}.
\end{equation}
This problem, including an anisotropic formulation introduced later, appears as a dimensional reduction in the free material design, see \cite{GRS}, and conductivity imaging, see \cite{JMN}. In this paper, we follow the approach to this problem from the point of view of geometric measure theory, following \cite{BGG}, \cite{JMN}, and \cite{SWZ}. In particular, we understand the boundary condition in the sense of traces of $BV$ functions.

This problem was introduced in \cite{SWZ}, where the authors estabilish that for continuous boundary data, under a set of conditions on an open bounded set $\Omega \subset \mathbb{R}^N$ slightly weaker than strict convexity, a unique solution to Problem (\ref{problem}) exists and it is continuous up to the boundary. However, if we relax some of these conditions, there arise additional problems:

(1) The first possible difficulty concerns discontinuous boundary data. Let us recall two results valid for $\Omega \subset \mathbb{R}^2$: as the example from \cite{ST} shows, if the boundary data are discontinuous on a set of positive measure, there may be no minimisers to Problem (\ref{problem}). On the other hand, as proved in \cite{Gor1}, for boundary data $f \in BV(\partial\Omega)$ there exists a minimiser to Problem (\ref{problem}); notice that in this case the set of discontinuities is countable. We will address this issue in Section \ref{sec:discontinuous}.

(2) The second possible difficulty concerns unbounded sets $\Omega \subset \mathbb{R}^N$. Even if the set $\Omega$ is strictly convex, the construction from \cite{SWZ} fails, as it involves minimalisation of perimeter in the class of sets which need not admit even one set with finite perimeter. Therefore, we need to work with approximations to both the set $\Omega$ and the boundary data $f$ in order to prove existence of minimisers. We will address this issue in Section \ref{sec:unbounded}.

(3) Finally, if we relax the assumptions concerning strict convexity of $\Omega$, the situation becomes much different - if $\Omega$ was convex with a flat part on the boundary, there exist continuous boundary data, for which there is no minimiser to Problem (\ref{problem}). Then, we need a different approach, involving finding a set of admissibility conditions sufficient for existence of minimisers. This issue is outside the scope of this paper and is explored for instance in \cite{RS}.

The purpose of this manuscript is twofold: firstly, we want to prove existence of minimisers for boundary data that may be discontinuous on a set of measure zero; for instance, when the boundary data $f \in BV(\partial\Omega)$ and the jump set of $f$ is small enough. Furthermore, we explore this problem in an anisotropic setting, where we encounter additional difficulties concerning our regularity and convexity assumptions on $\phi$ and $\partial\Omega$ due to the use of the barrier condition and comparison principles.

Secondly, we use the technique developed in the first part to find the appropriate function space on $\partial\Omega$, for which the least gradient problem (and its anisotropic versions) are well-posed for unbounded domains. For geometrical reasons arising even in the bounded domain case, we assume $\Omega$ to be strictly convex and not equal to $\mathbb{R}^N$. We want to deal with two kinds of phenomena: the regularity of boundary data and the shape of the domain. The main existence result, Theorem \ref{thm:generalexistence}, works in quite general setting, but imposing additional constraints on both the domain and boundary data gives us uniqueness and additional regularity of minimisers. \\

Finally, in both cases we will additionally address the anisotropic case. We are interested in the following version of the least gradient problem:
\begin{equation}\label{aproblem}\tag{ALGP}
\min \{ \int_\Omega \phi(x,Du), \quad u \in BV(\Omega), \quad u|_{\partial\Omega} = f  \},
\end{equation}
where $\phi$ will be a 1-homogenous function convex in the second variable, which either depends only on the second variable (in other words, it is a norm on $\mathbb{R}^N$) or satisfies some regularity assumptions implying a comparison principle. In particular, the second case covers the weighted least gradient problem, where $\phi(x,Du) = a(x)|Du|$ for $a \in C^{1,1}(\Omega)$. These two cases will require slightly different methods; in the former, we will exploit the translation invariance and projections, while in the latter we will rely primarily on the comparison principle.

\section{Preliminaries}

\subsection{Least gradient functions}

In this section, we recall the definition of least gradient functions and their basic properties.

\begin{dfn}\label{dfn:lg}
Let $\Omega \subset \mathbb{R}^N$ be open. We say that $u \in BV(\Omega)$ is a function of least gradient, if for every $v \in BV(\Omega)$ compactly supported in $\Omega$ we have
\begin{equation*}
\int_\Omega |Du| \leq \int_\Omega |D(u + v)|.
\end{equation*}
In case when $\Omega$ is bounded with Lipschitz boundary, this is equivalent to the condition that $v \in BV_0(\Omega)$, see \cite[Theorem 2.2]{SZ0}. This equivalence is proved using an approximation with functions of the form $v_n = v \chi_{\Omega_n}$ for suitably chosen $\Omega_n$ and the proof does not extend well to the case when $\Omega$ is unbounded.
\end{dfn}

\begin{dfn}
We say that $u \in BV_{loc}(\Omega)$ is a solution to Problem (\ref{problem}), if $u$ is a function of least gradient and the trace of $u$ equals $f$, i.e. for almost every $x \in \partial\Omega$ we have 
$$ \dashint_{B(x,r) \cap \Omega} |f(x) - u(y)| dy = 0.$$
We prefer to state the trace condition in this way in order to avoid discussion on existence and continuity of the trace operator for unbounded sets $\Omega$.
\end{dfn}

Now, we recall three classical theorems on least gradient functions. The first one is Miranda's theorem on stability of least gradient functions:

\begin{tw}\label{thm:miranda}
$($\cite[Theorem 3]{Mir}$)$
Let $\Omega \subset \mathbb{R}^N$ be open. Suppose $\{ u_n \} \subset BV(\Omega)$ is a sequence of least gradient functions in $\Omega$ convergent in $L^1_{loc}(\Omega)$ to $u \in BV(\Omega)$. Then $u$ is a function of least gradient in $\Omega$. \qed
\end{tw}

However, this theorem has a very important limitation: even if $\Omega$ is bounded, as the trace operator is not continuous with respect to $L^1$ convergence, it gives us no control on the trace of the limit function $u$. In fact, in this manuscript and other works concerning least gradient functions (\cite{GRS}, \cite{SWZ}) much effort is devoted to prove that the trace of the limit function is correct.

The second result is a theorem by Bombieri-de Giorgi-Giusti, which gives us a link between the function $u$ of least gradient and the regularity of its superlevel sets. Here and in the whole manuscript let us denote $E_t = \{ u \geq t \}$.

\begin{tw}\label{tw:bgg}
(\cite[Theorem 1]{BGG}) Suppose that $\Omega \subset \mathbb{R}^N$ is open and let $u \in BV(\Omega)$ be a function of least gradient in $\Omega$. Then for every $t \in \mathbb{R}$ the set $E_t$ is minimal in $\Omega$, i.e. the function $\chi_{E_t}$ is of least gradient.
\end{tw}

The third result, by Sternberg-Williams-Ziemer, concerns the existence of minimisers to Problem (\ref{problem}) in case when $\Omega \subset \mathbb{R}^N$ is an open bounded strictly convex set. While the assumptions on $\Omega$ are a little weaker (the authors assume that $\Omega$ has positive weak mean curvature on a dense set and that $\partial\Omega$ is not locally area-minimising), for simplicity we are interested only in strictly convex sets.

\begin{tw}\label{thm:swz}
(\cite[Theorem 4.1]{SWZ}) Let $\Omega \subset \mathbb{R}^N$ be an open, bounded, strictly convex set and suppose that $f\in C(\partial\Omega)$. Then there exists a unique minimiser $u \in BV(\Omega)$ to Problem (\ref{problem}) and additionally $u \in C(\overline{\Omega})$.
\end{tw}

Finally, as least gradient functions are $BV$ functions, they are defined up to a set of measure zero, we have to choose a proper representative if we want to state any regularity results. In this paper, following \cite{SWZ}, we employ the convention that a set of a bounded perimeter consists of all its points of positive density.

\subsection{Anisotropic formulation}

Firstly, we recall the notion of a metric integrand and $BV$ spaces with respect to a metric integrand. We follow the construction in \cite{AB}.

\begin{dfn}
Let $\Omega \subset \mathbb{R}^N$ be an open bounded set with Lipschitz boundary. A continuous function $\phi: \overline{\Omega} \times \mathbb{R}^N \rightarrow [0, \infty)$ is called a metric integrand, if it satisfies the following conditions: \\
\\
$(1)$ $\phi$ is convex with respect to the second variable for a.e. $x \in \overline{\Omega}$; \\
$(2)$ $\phi$ is homogeneous with respect to the second variable, i.e.
\begin{equation*}
\forall \, x \in \overline{\Omega}, \quad \forall \, \xi \in \mathbb{R}^N, \quad \forall \, t \in \mathbb{R} \quad \phi(x, t \xi) = |t| \phi(x, \xi);
\end{equation*}
$(3)$ $\phi$ is bounded and elliptic in $\overline{\Omega}$, i.e.
\begin{equation*}
\exists \,  \lambda, \Lambda  > 0 \quad \forall \, x \in \overline{\Omega}, \quad \forall \, \xi \in \mathbb{R}^N \quad \lambda |\xi| \leq \phi(x, \xi) \leq \Lambda |\xi|.
\end{equation*}
These conditions apply to most cases considered in the literature, such as the classical least gradient problem, i.e. $\phi(x, \xi) = |\xi|$ (see \cite{SWZ}), the weighted least gradient problem, i.e. $\phi(x, \xi) = g(x) |\xi|$ (see \cite{JMN}), where $g \geq c > 0$, and $l_p$ norms for $p \in [1, \infty]$, i.e. $\phi(x, \xi) = \| \xi \|_p$ (see \cite{Gor1}).
\end{dfn}

\begin{dfn}
The polar function of $\phi$ is $\phi^0: \overline{\Omega} \times \mathbb{R}^N \rightarrow [0, \infty)$ defined as

\begin{equation*}
\phi^0 (x, \xi^*) = \sup \, \{ \langle \xi^*, \xi \rangle : \, \xi \in \mathbb{R}^N, \, \phi(x, \xi) \leq 1 \}.
\end{equation*} 
\end{dfn}

\begin{dfn}
Let $\phi$ be a metric integrand continuous and elliptic in $\overline{\Omega}$. For a given function $u \in L^1(\Omega)$ we define its $\phi-$total variation in $\Omega$ by the formula:

\begin{equation*}
\int_\Omega |Du|_\phi = \sup \, \{ \int_\Omega u \, \mathrm{div} \, \mathbf{z} \, dx : \, \phi^0(x,\mathbf{z}(x)) \leq 1 \, \, \, \text{a.e.}, \quad \mathbf{z} \in C_c^1(\Omega)  \}.
\end{equation*}
Another popular notation for the $\phi-$total variation is $\int_\Omega \phi(x, Du)$. We will say that $u \in BV_\phi(\Omega)$ if its $\phi-$total variation is finite; furthermore, let us define the $\phi-$perimeter of a set $E$ as
$$P_\phi(E, \Omega) = \int_{\Omega} |D\chi_E|_\phi.$$
If $P_\phi(E, \Omega) < \infty$, we say that $E$ is a set of bounded $\phi-$perimeter in $\Omega$.
\end{dfn}

\begin{uw}
By property $(3)$ of a metric integrand
$$\lambda \int_\Omega |Du| \leq \int_\Omega |Du|_\phi \leq \Lambda \int_\Omega |Du|.$$
In particular, $BV_\phi(\Omega) = BV(\Omega)$ as sets; however, they are equipped with different (but equivalent) norms and corresponding strict topologies. Moreover, the space $BV_{\phi}(\Omega)$ satisfies the same basic properties as the isotropic space $BV(\Omega)$, such as lower semicontinuity of the $\phi-$total variation with respect to convergence in $L^1$, the co-area formula, and the approximation by smooth functions in the strict topology:
\end{uw}

\begin{uw}\label{lem:scisaniz}
Suppose that $\Omega$ is an open bounded set with Lipschitz boundary and $\phi$ is a metric integrand continuous and elliptic in $\overline{\Omega}$. Let $v \in BV_\phi(\Omega)$ and $Tv = f$. Then there exists a sequence $v_n \in C^\infty(\Omega) \cap BV(\Omega)$ such that $v_n \rightarrow v$ strictly in $BV_\phi(\Omega)$ and $Tv_n = f$ (in the isotropic case, see \cite[Corollaries 1.17, 2.10]{Giu}). \qed
\end{uw}

Finally, we will use the following integral representation of the $\phi-$total variation (\cite{AB}, \cite{JMN}):

\begin{stw}\label{stw:repcalkowa}
Let $\varphi: \overline{\Omega} \times \mathbb{R}^N \rightarrow \mathbb{R}$ be a metric integrand. Then we have an integral representation:

\begin{equation*}
\int_\Omega |Du|_\phi = \int_\Omega \phi(x, \nu^u(x)) \, |Du|,
\end{equation*}
where $\nu^u$ is the Radon-Nikodym derivative $\nu^u = \frac{d Du}{d |Du|}$. If we take $u$ to be a characteristic function of a set $E$ with a $C^1$ boundary, we have

\begin{equation*}
P_\phi(E, \Omega) = \int_\Omega \phi(x, \nu_E) \, d \mathcal{H}^{N-1},
\end{equation*}
where $\nu(x)$ is the (Euclidean) unit vector normal to $\partial E$ at $x \in \partial E$. Let us also denote by $\tau(x)$ the unit vector tangent to $\partial E$ at $x \in \partial E$. \qed
\end{stw}

\subsection{$\phi-$least gradient functions}

Now, we turn our attention to the precise formulation of Problem (\ref{aproblem}). Then we recall several known properties of the minimisers and discuss how will the geometry of $\Omega$ and regularity of $\phi$ come into play.

\begin{dfn}
Let $\Omega \subset \mathbb{R}^N$ be an open bounded set with Lipschitz boundary. We say that $u \in BV_\phi(\Omega)$ is a function of $\phi-$least gradient, if for every compactly supported $v \in BV_\phi(\Omega)$ we have

\begin{equation*}
\int_\Omega |Du|_\phi \leq \int_\Omega |D(u + v)|_\phi.
\end{equation*}
We say that $u$ is a solution to Problem (\ref{aproblem}), the anisotropic least gradient problem with boundary data $f$, if $u$ is a function of $\phi-$least gradient and $Tu = f$.
\end{dfn}

We will recall a few properties of functions of $\phi-$least gradient. Firstly, we state an anisotropic version of Theorem \ref{tw:bgg}. Its proof in both directions is based on the the co-area formula.

\begin{stw}\label{stw:anizobgg}
(\cite[Theorem 3.19]{Maz}) Let $\Omega \subset \mathbb{R}^N$ be an open bounded set with Lipschitz boundary. Assume that the metric integrand $\phi$ has a continuous extension to $\mathbb{R}^N$. Take $u \in BV_\phi(\Omega)$. Then $u$ is a function of $\phi-$least gradient in $\Omega$ if and only if $\chi_{\{ u > t \}}$ is a function of $\phi-$least gradient for almost all $t \in \mathbb{R}$. \qed
\end{stw}

Both in the isotropic and anisotropic case existence and uniqueness of minimisers depend on the geometry of $\Omega$. Suppose that the boundary data are continuous. In the isotropic case, the necessary and sufficient condition was introduced in \cite{SWZ} and in two dimensions it is equivalent to strict convexity of $\Omega$. In the anisotropic case, a sufficient condition (see \cite[Theorem 1.1]{JMN}) is the {\it barrier condition}:

\begin{dfn}\label{def:barrier}
(\cite[Definition 3]{JMN}) Let $\Omega \subset \mathbb{R}^N$ be an open bounded set with Lipschitz boundary. Suppose that $\phi$ is an elliptic metric integrand. We say that $\Omega$ satisfies the barrier condition (with respect to $\phi$) if for every $x_0 \in \partial\Omega$ and sufficiently small $\varepsilon > 0$, if $V$ minimises $P_\phi(\cdotp; \mathbb{R}^N)$ in 

$$ \{ W \subset \Omega: W \backslash B(x_0, \varepsilon) = \Omega \backslash B(x_0, \varepsilon) \} $$
then
$$\partial V \cap \partial\Omega \cap B(x_0, \varepsilon) = \emptyset.$$
In the isotropic case $\phi(x, \xi) = \| \xi \|_2$ this is equivalent, at least for sets with $C^2$ boundary, to the condition introduced in \cite{SWZ}.
\end{dfn}

\begin{tw}\label{jmn:existence}
(\cite[Theorem 1.1]{JMN}) Suppose that $\phi$ is a metric integrand and $\Omega \subset \mathbb{R}^N$ be an open bounded set with Lipschitz boundary which satisfies the barrier condition. Let $f \in C(\partial\Omega)$. Then there exists a minimiser $u \in BV_\phi(\Omega)$ to Problem (\ref{aproblem}).
\end{tw}

Our second concern in the study of least gradient functions is the regularity of $\phi-$minimal sets; it is related to the maximum and comparison principles for $\phi-$minimal sets. Explicitly, we will assume that
\begin{equation*}\tag{H}\label{hypothesis}
\text{If } E \text{ is a }\phi-\text{minimal set, then } \twopartdef{\mathcal{H}^{n-3}( \text{sing } \partial E) < \infty}{n \geq 4}{\text{sing } \partial E = \emptyset}{n \leq 3}
\end{equation*}
i.e. $\partial E$ is of class $C^2$ apart from a set of finite $\mathcal{H}^{n-3}$ measure. In dimensions up to three, the singular set is in fact empty. In fact, it is a regularity assumption on $\phi$; the conditions on $\phi$ which imply $(H)$ involve uniform convexity and regularity somewhat weaker than $C^1$ in the first (spatial) variable and $C^3$ in the second (directional) variable. The sufficiency of these conditions is proved in \cite{SSA} (a detailed discussion can be found in \cite{JMN}, where the authors additionally show that one cannot relax the regularity in the spatial variable). 

\begin{tw}\label{jmn:comparison}
(\cite[Theorem 1.2]{JMN}) Suppose that $\Omega \subset \mathbb{R}^N$ is an open bounded set with connected Lipschitz boundary. Let $\phi$ be a metric integrand which additionally satisfies hypothesis (\ref{hypothesis}). Let $f_1, f_2 \in C(\partial\Omega)$ such that $f_1 \geq f_2$ on $\partial\Omega$. Let $u_i$ be minimisers to Problem (\ref{aproblem}) corresponding to $f_i$. Then $u_1 \geq u_2$ in $\Omega$. In particular, minimisers to Problem (\ref{aproblem}) with continuous boundary data are unique.
\end{tw}

\subsection{Technical lemmas}

Now, we write a very simple observation:

\begin{lem}\label{lem:variationestimate}
Let $\Omega \subset \mathbb{R}^N$ be an open bounded set with Lipschitz boundary. Let $u \in BV_\phi(\Omega)$ be a function of $\phi-$least gradient in $\Omega$. Then
$$ \int_\Omega |Du|_\phi \leq \int_{\partial\Omega} \phi(x, \nu) |Tu| \, d\mathcal{H}^{N-1}.$$
\end{lem}

\begin{dd}
As shown in \cite{Maz}, the relaxed functional for the anisotropic least gradient problem with boundary data $f$ is
$$ F_\phi(v) = \int_{\Omega} |Dv|_\phi + \int_{\partial\Omega} \phi(x, \nu) |Tv - f|.$$
As $u$ is a function of $\phi-$least gradient with boundary data $Tu$, we have
$$ \int_{\Omega} |Du|_\phi = \int_{\Omega} |Du|_\phi + \int_{\partial\Omega} \phi(x, \nu) |Tu - Tu| d\mathcal{H}^{N-1} = F_\phi(u) \leq $$
$$ \leq F_\phi(v \equiv 0) = \int_{\Omega} 0 \,dx + \int_{\partial\Omega} \phi(x, \nu) |0 - Tu| \, d\mathcal{H}^{N-1} = \int_{\partial\Omega} \phi(x, \nu) |Tu| \, d\mathcal{H}^{N-1}.$$
\qed
\end{dd}

And another simple observation, based on the Poincar\'{e} inequality:

\begin{lem}\label{lem:poincare}
Let $\Omega$ be an open bounded set with Lipschitz boundary which lies in a strip of width $d$ and let $u \in BV(\Omega)$. Then
$$\| u \|_{L^1(\Omega)} \leq C(d) (\int_{\Omega} |Du| + \int_{\partial\Omega} |Tu| d\mathcal{H}^{N-1}).$$
\end{lem}

\begin{dd}
Let us extend the function $u$ by $0$ on $\mathbb{R}^N \backslash \Omega$. The extension $\widetilde{u}$ is a function with compact support in $\mathbb{R}^N$ and by the Poincar\'{e} inequality for $\widetilde{u}$ we have
$$ \| u \|_{L^1(\Omega)} = \| \widetilde{u} \|_{L^1(\mathbb{R}^N)} \leq C(d) \int_{\mathbb{R}^N} |D\widetilde{u}| = C(d) (\int_\Omega |Du| + \int_{\partial\Omega} |Tu| \, d\mathcal{H}^{N-1}) + \int_{\mathbb{R}^N \backslash \overline{\Omega}} 0,$$
where $C(d)$ is the constant in the Poincar\'{e} inequality which depends only on the width $d$ of the strip containing $\Omega$. \qed
\end{dd}

\begin{lem}\label{lem:gor1conv}
Let $\Omega \subset \mathbb{R}^N$ be an open bounded set with Lipschitz boundary. Suppose that $\phi$ is a metric integrand and let $\{ f_n \} \subset L^1(\partial\Omega)$ be a uniformly bounded sequence. Finally, let $u_n \in BV_\phi(\Omega)$ be functions of $\phi-$least gradient with traces $f_n$. Then $u_n$ has a convergent subsequence $u_{n_k} \rightarrow u$ in $L^1(\Omega)$.
\end{lem}
The above Lemma generalises \cite[Proposition 4.1]{Gor1} concerning the isotropic case. Furthermore, if $f_n \rightarrow f$ in $L^1(\partial\Omega)$, it does not imply that $Tu = f$, as the trace operator is not continuous with respect to convergence in $L^1(\Omega)$.

\begin{dd}
Recall that $BV_\phi(\Omega) = BV(\Omega)$ as sets. We start by estimating the $L^1$ norm of $u$ using the Poincar\'{e} inequality. By Lemma \ref{lem:poincare} we have
$$ \| u_n \|_{L^1(\Omega)} \leq C(\Omega) (\int_\Omega |Du_n| + \int_{\partial\Omega} |Tu_n| \, d\mathcal{H}^{N-1}).$$
We estimate the (isotropic) total variations of $u_n$ using Lemma \ref{lem:variationestimate} and the equivalence of norms between $BV_\phi(\Omega)$ and $BV(\Omega)$:
$$ \int_\Omega |Du_n| \leq \lambda^{-1} \int_\Omega |Du_n|_\phi \leq \lambda^{-1} \int_{\partial\Omega} \phi(x, \nu) |Tu_n| \, d\mathcal{H}^{N-1} \leq \lambda^{-1} \Lambda \int_{\partial\Omega} |Tu_n| \, d\mathcal{H}^{N-1}.$$
We bring these two estimates together and get
$$ \| u_n \|_{BV(\Omega)} = \| u_n \|_{L^1(\Omega)} + \int_\Omega |Du_n| \leq (C(\Omega) + 1) (\int_\Omega |Du_n| + \int_{\partial\Omega} |Tu_n| \, d\mathcal{H}^{N-1}) \leq$$
$$ \leq (C(\Omega) + 1) (\lambda^{-1} \Lambda + 1) \int_{\partial\Omega} |Tu_n| d\mathcal{H}^{N-1} = (C(\Omega) + 1) (\lambda^{-1} \Lambda + 1) \int_{\partial\Omega} |f_n| d\mathcal{H}^{N-1}.$$
As the sequence $f_n$ is uniformly bounded in $L^1(\partial\Omega)$, the sequence $u_n$ is uniformly bounded in $BV(\Omega)$, hence it admits a convergent subsequence in $L^1(\Omega)$. \qed
\end{dd}

\begin{lem}\label{lem:swzbdry}
(\cite[Lemma 3.3]{SWZ}) Suppose that $\Omega \subset \mathbb{R}^N$ is an open, bounded, strictly convex set and let $g \in C(\partial\Omega)$ be continuous. Let $u \in BV(\Omega)$ be a least gradient function with trace $g$ provided by Theorem \ref{thm:swz}. Then for almost all $t \in [a,b]$ we have $\partial E_t \cap \partial\Omega \subset g^{-1}(t)$.
\end{lem}

\section{Results for discontinuous boundary data}\label{sec:discontinuous}

This section is devoted to proving existence of minimisers for bounded sets, but with fairly general boundary data. We will consider two different cases. Firstly, we will consider the case when $\phi$ is a norm with strictly convex unit ball, without any regularity assumptions on $\phi$. Secondly, we will assume that the metric integrand $\phi$ may depend on location, but has to satisfy the regularity hypothesis (\ref{hypothesis}). In particular, both approaches cover the isotropic least gradient problem.

\subsection{Existence theorems}

The main results of this Section, Theorem \ref{thm:discexistencev1} and Theorem \ref{thm:discexistencev2}, concern the case when the discontinuity set of the boundary data $\phi$ is small - the precise assumption is that its $\mathcal{H}^{N-1}-$measure is zero. These results are motivated by \cite[Theorem 1.1]{Gor1}, which states that in two dimensions, if $\Omega$ is strictly convex and $f \in BV(\partial\Omega)$, then there exists a minimiser to Problem (\ref{problem}). The two Theorems extend this result to higher dimensions, while generalizing it to anisotropic cases as well.

\begin{tw}\label{thm:discexistencev1}
Let $\phi$ be a norm on $\mathbb{R}^N$ such that the unit ball $B_\phi(0,1)$ is strictly convex and let $\Omega \subset \mathbb{R}^N$ be a strictly convex set which satisfies the barrier condition with respect to $\phi$. Suppose that $L^1(\partial\Omega)$ is a function such that $\mathcal{H}^{N-1}$-almost all points of $\partial\Omega$ are continuity points of $f$. Then there exists a minimiser to the Problem (\ref{problem}) with boundary data $f$.
\end{tw}

\begin{dd}
1. We want to define a sequence of approximations $f_n$ which is continuous, converges almost uniformly to $f$ and which locally preserves the $L^\infty$ bounds of $f$. Mollification has all of the above properties, but it need not be defined if $\partial\Omega$ is not a Lie group; therefore we have to construct a similar operator.

Let $\rho_\varepsilon$ be a standard mollification kernel on $\mathbb{R}^{N-1}$ with the diameter of support equal to $\varepsilon$. As $\Omega$ is strictly convex and bounded, its boundary is compact and locally it is a graph of a Lipschitz function. Hence, it is a topological manifold, equipped with an atlas of finitely many (due to compactness) bi-Lipschitz maps $\phi_m: U_n \rightarrow \partial\Omega$, where $U_m \subset \mathbb{R}^{N-1}$ is open (the inverse of each $\phi_m$ is a projection, so also a Lipschitz map). We denote $V_m = \phi_m(U_m)$; the sets $V_n$ form an open cover of $\partial\Omega$. Let $\varphi_m$ be a continuous partition of unity subject to the cover $V_m$. 

We define $f_n \in C(\partial\Omega)$ in the following way: in the domain of a map $\phi_i$, we pull back $f$ with a map $\phi_i$ to $\mathbb{R}^{N-1}$. We mollify the pullback with a kernel $\rho_{\frac{1}{n}}$ and go back to $\partial\Omega$. In other words, we write
$$ f_n(x) = \sum_{i = 1}^m \varphi_i(x) g_i(\phi_i^{-1}(x)), $$
where
$$ g_i = (\varphi_i \circ \phi_i)(f \circ \phi_i) * \rho_{\frac{1}{n}} \in C_c^\infty(U_i),$$
and whenever $\phi_i^{-1}(x)$ is not defined, then $\varphi_i$ equals zero and nothing changes in the sum. Then $f_n$ is a continuous function and the value of $f_n$ at $x$ depends only on the values of $f$ in a ball $B(x, r(n))$, where $r(n) \rightarrow 0$ as $n \rightarrow \infty$; the precise form of $r(n)$ depends on the Lipschitz constants of maps $\phi_m$.

2. As $f_n \in C(\partial\Omega)$, there exist solutions $u_n \in BV(\Omega) \cap C(\partial\Omega)$ of the least gradient problem. By Lemma \ref{lem:gor1conv} they converge on a subsequence to a function $u$ of least gradient in $\Omega$; we only have to ensure that the trace is correct.

3. The set of points on $\partial\Omega$ where the trace of $u$ is well-defined by the mean value property is of $\mathcal{H}^{N-1}-$full measure. Furthermore, by our assumption $\mathcal{H}^{N-1}-$almost all points of $\partial\Omega$ are continuity points of $f$. We denote the set (of $\mathcal{H}^{N-1}-$full measure) where the two above points hold by $\mathcal{Z} \subset \partial\Omega$.

4. Fix $x_0 \in \mathcal{Z}$; in particular $x_0$ is a point of continuity of $f$. Fix any $\delta > 0$. Then there exists a ball $B(x_0,r)$ such that for all $x \in B(x_0,r) \cap \partial\Omega$ we have
$$ f(x_0) - \delta \leq f(x) \leq f(x_0) + \delta.$$
As the construction of the approximating sequence $f_n$ involves mollification, we have for sufficiently large $n$ (so that the mollification kernel has support with diameter such that $r(n) < \frac{r}{2}$)
$$ f(x_0) - \delta \leq f_n(x) \leq f(x_0) + \delta$$
on a smaller set $B(x_0, \frac{r}{2}) \cap \partial\Omega$.

5. We introduce the following notation: let $H$ be a supporting hyperplane at $x_0$. We choose it from the set of all supporting hyperplanes so that the normal $\nu$ to $H$ at $x_0$ points inside $\Omega$. Let $H^+$ be the halfspace bounded by $H$ which does not contain $\Omega$. Take $s > 0$ small enough, so that $(H + s\nu) \cap \Omega \subset \subset B(x_0, \frac{r}{2})$. We want to estimate the value of $u_n$ in the set $\Omega' = (H_+ + s\nu) \cap \Omega$.

6. We will see that in the set $\Omega'$ for sufficiently large $n$ all the functions $u_n$ satisfy the inequality
$$ f(x_0) - \delta \leq u_n(x) \leq f(x_0) + \delta.$$
Suppose otherwise; without loss of generality let $u_n(x') = t > f(x_0) + \delta$. By the continuity of $u$ the set $\partial E_t$ intersects $\Omega'$. Furthermore, again by the continuity of $u$ the set $\partial E_s$ intersects $\Omega'$ on some interval $\mathcal{I}$ close to $t$.

7. Fix some $s \in \mathcal{I}$ for which statement of Lemma \ref{lem:swzbdry} holds (it is possible, as it holds for almost all $t \in \mathbb{R}$). Then the set $E_s$ is closed in $\Omega$, intersects $\Omega'$, but by Lemma \ref{lem:swzbdry} it does not intersect $\partial \Omega' \cap \partial \Omega$. Hence $E_s$ has positive distance to the closed set $\partial \Omega' \cap \partial \Omega$ and the function $\chi_{E_s \backslash \Omega'}$ has the same trace as the function $\chi_{E_s}$; however, the former function has strictly smaller $\phi-$perimeter, as any parts of $E_s$ intersecting $\Omega'$ are projected onto the hyperplane $H + s \nu$. This contradicts Theorem \ref{tw:bgg}, hence $u_n(x) \in [f(x_0) - \delta, f(x_0) + \delta]$.

8. Now, we see that the trace of $u$ at $x_0$ equals $f$. Passing to the pointwise limit almost everywhere in the inequality from the Step 6 (possibly passing to a subsequence), we obtain that in a ball $B(x_0, \rho) \subset \Omega'$ we have
$$ f(x_0) - \delta \leq u(x) \leq f(x_0) + \delta.$$
As for arbitrary $\delta > 0$ there exists a ball $B(x_0, \rho)$ for which the above inequality is satisfied, we see that
$$ \lim_{\rho \rightarrow 0}\, \,  \text{ess sup}_{B(x_0,\rho)} \, |u(y) - f(x_0)| = 0,$$
so $Tu(x_0) = f(x_0)$. This equality holds for $\mathcal{H}^{N-1}$-almost all $x_0 \in \partial\Omega$, so $Tu = f$. \qed
\end{dd}

Now, in the following remarks we will discuss the necessity of the assumptions in the above Theorem. There are two types of assumptions in play: the barrier condition and its relation to strict convexity; and the fact that $B_\phi(0,1)$ is strictly convex.

\begin{uw}
The assumption that the set is both strictly convex and satisfies the barrier condition with respect to $\phi$ is quite natural for two reasons:

1. Firstly, as the hyperplane is a minimal surface for any norm $\phi$, if a set satisfies a barrier condition for a norm $\phi$, it satisfies the barrier condition for the isotropic norm $l^2$; but the barrier condition for the isotropic norm is equivalent to the fact that the boundary of $\Omega$ has positive mean curvature on a dense subset and is not locally-area minimising, which is something stronger than convexity and a little weaker than strict convexity.

2. Secondly, let us highlight what does not work when the set is not strictly convex. The only problem, maybe only technical, is choosing the appropriate hyperplane - so that the set $\Omega'$ is indeed contained in some ball; for instance, in three dimensions, if $\partial\Omega$ contains a line segment, then the set $\Omega'$ would contain a neighbourhood of this line segment and not be contained in any $B(x_0,r)$ for small $r$.
\end{uw}

\begin{uw}
If $\partial \Omega$ is a smooth manifold with group structure, the first part of the proof becomes simpler - we can just take the approximating sequence to be defined using convolution with a mollification kernel: $f_n = f * \rho_{\frac{1}{n}}$. Moreover, in Step 5 we take simply a tangent hyperplane instead of choosing a proper supporting hyperplane.
\end{uw}

\begin{uw}
Let us highlight what does and does not work when $B_\phi(0,1)$ is not strictly convex. Then, as we can see in \cite{Gor3}, no set with $C^1$ boundary satisfies the barrier condition; thus we do not know if there exist minimisers for continuous boundary data. Furthermore, the projection on a hyperplane does not necessarily decrease the anisotropic total variation. We will further discuss the two-dimensional case in Proposition \ref{stw:existence2d}.
\end{uw}

\begin{uw}
A similar method has been developed in dimension two in order to prove existence of minimisers in two situations: in the isotropic least gradient problem, when the domain is not strictly convex, but the boundary data satisfy some admissibility conditions, see \cite{RS}; and when the domain is strictly convex, but the norm $\phi$ is such that $B_\phi(0,1)$ is not strictly convex, see \cite{Gor3}.
\end{uw}

Now, we turn our attention to the case when $\phi$ is a metric integrand and not necessarily a norm; however, we will assume high regularity of $\phi$ in order to use a comparison principle. We begin by a lemma which generalises a very simple observation from the isotropic case: 

\begin{lem}\label{lem:uniformbarrier}
Let $\phi$ be a metric integrand satisfying hypothesis (\ref{hypothesis}). Suppose that $\Omega \subset \mathbb{R}^N$ is an open set with $C^2$ boundary which satisfies the barrier condition. Let $g \in C(\partial\Omega)$. If $g$ is constant in some neighbourhood of $x_0$ in $\partial\Omega$, then the minimiser $v$ with boundary data $g$ is constant in some neighbourhood of $x_0$ in $\Omega$.
\end{lem}

\begin{dd}
Firstly, let us notice that if $V_\alpha$ is a family of sets from the definition of barrier condition, i.e. for each $\alpha \in \mathcal{A}$ the set $V_\alpha$ minimises $P_\phi(\cdot, \mathbb{R}^N)$ in
\begin{equation*}\label{eq:class}\tag{*}
\{ W \subset \Omega: \quad W \backslash B(x_0, \varepsilon) = \Omega \backslash B(x_0, \varepsilon)\}.
\end{equation*}
Consider the set $V = \bigcup_{\alpha \in \mathcal{A}}$. It satisfies $V \backslash B(x_0, \varepsilon) = \Omega \backslash B(x_0, \varepsilon)$ and by lower semicontinuity of the $\phi-$total variation it minimises the perimeter in the class (\ref{eq:class}). By the barrier condition $\partial V \cap \partial \Omega \cap B(x_0, \varepsilon) = \emptyset$; but as $\partial\Omega \in C^2$ and $\phi$ satisfies (\ref{hypothesis}), the distance from $x_0$ to $V$ is positive, hence $V \cap B(x_0,r) = \emptyset$ for some $r = r(\varepsilon) > 0$. 

Now, suppose that $g$ is constant in $B(x_0, \varepsilon) \cap \partial\Omega$ with constant value $t_0$. Then take $V_t = \{ v > t \}$ with $t > t_0$. The set $V_t$ falls into the class (\ref{eq:class}) and by the continuity of $g$ its closure does not contain $x_0$. Hence, there is a positive distance $r = r(\varepsilon)$ separating $x_0$ and $V_t$ for all $t > t_0$, so $v \leq t_0$ in $B(x_0, r(\varepsilon)) \cap \Omega$. Using a similar argument with $\widetilde{V}_t = \{ v < t \}$ with $t < t_0$ we obtain that $v \geq t_0$ in $B(x_0, \varepsilon) \cap \Omega$. Hence $v$ is constant with value $t_0$ in some neighbourhood of $x_0$ in $\Omega$. \qed
\end{dd}

\begin{tw}\label{thm:discexistencev2}
Let $\phi$ be a metric integrand satisfying hypothesis (\ref{hypothesis}) and let $\Omega \subset \mathbb{R}^N$ be an open bounded set with $C^2$ boundary which satisfies the barrier condition. Suppose that $L^1(\partial\Omega)$ is a function such that $\mathcal{H}^{N-1}$-almost all points of $\partial\Omega$ are continuity points of $f$. Then there exists a minimiser to the Problem (\ref{aproblem}) with boundary data $f$.
\end{tw}

\begin{dd}
1. We define $f_n = f * \rho_{\frac{1}{n}}$, where $\rho_\varepsilon$ is a standard mollification kernel with the diameter of support equal to $\varepsilon$. Then $f_n$ is a smooth function and the value of $f_n$ at $x$ depends only on the values of $f$ in the ball $B(x, \frac{1}{n})$.

2. As $f_n \in C^\infty(\partial\Omega)$ and $\Omega$ satisfies the barrier condition, there exist solutions $u_n \in BV(\Omega) \cap C(\Omega)$ of the least gradient problem. By Lemma \ref{lem:gor1conv} they converge on a subsequence to a function $u$ of $\phi-$least gradient in $\Omega$; we only have to ensure that the trace is correct.

3. The set of points on $\partial\Omega$ where the trace of $u$ is well-defined by the mean value property is of $\mathcal{H}^{N-1}-$full measure. Furthermore, by our assumption $\mathcal{H}^{N-1}-$almost all points of $\partial\Omega$ are continuity points of $f$. We denote the set (of $\mathcal{H}^{N-1}-$full measure) where the two above points hold by $\mathcal{Z} \subset \partial\Omega$.

4. Fix $x_0 \in \mathcal{Z}$; in particular $x_0$ is a point of continuity of $f$. Fix any $\delta > 0$. Then there exists a ball $B(x_0,r)$ such that for all $x \in B(x_0,r) \cap \partial\Omega$ we have
$$ f(x_0) - \delta \leq f(x) \leq f(x_0) + \delta.$$
As the sequence $f_n$ is constructed via mollification, for sufficiently large $n$ (so that the mollification kernel has support with diameter smaller than $\frac{r}{2}$)
$$ f(x_0) - \delta \leq f_n(x) \leq f(x_0) + \delta$$
on a smaller set $B(x_0, \frac{r}{2}) \cap \partial\Omega$.

5. We will use the comparison principle to obtain a uniform bound on the sequence $u_n$ in some neighbourhood of $x_0$. To this end, let us define $f_n^\pm \in C(\partial\Omega)$ as follows:
$$ f_n^+ = \max(f_n, f(x_0) + \delta), \qquad f_n^-=\min(f_n, f(x_0) + \delta).$$
By definition $f_n^- \leq f_n \leq f_n^+$ on $\partial\Omega$ and $u_n^\pm$ are constant in $B(x_0, \frac{r}{2}) \cap \partial\Omega$ and equal to $f(x_0) \pm \delta$ respectively. As $\Omega$ satisfies the barrier condition, by Theorem \ref{jmn:existence} there exist minimisers $u_n^\pm$ for boundary data $f_n^\pm$; as $\phi$ satisfies hypothesis (\ref{hypothesis}), by Theorem \ref{jmn:comparison} we have
$$ u_n^- \leq u_n \leq u_n^+ \qquad \text{ in } \Omega.$$

6. By Lemma \ref{lem:uniformbarrier} the barrier condition implies that if boundary data $g$ is constant in some neighbourhood of $x_0$ in $\partial\Omega$, then the minimiser $v$ with boundary data $g$ is constant in some neighbourhood of $x_0$ in $\Omega$. We apply this to the sequence $u_n$ and obtain a ball $B(x_0,r')$ such that $u_n^\pm = f(x_0) \pm \delta$ in $\Omega \cap B(x_0,r')$. We notice that the radius $r'$ is independent of $n$ and depends only on the original radius $r$.

7. By the comparison principle (Theorem \ref{jmn:comparison}) all the functions $u_n$ satisfy the inequality
$$ f(x_0) - \delta \leq u_n(x) \leq f(x_0) + \delta$$
in $B(x_0, r')$. Now, we see that the trace of $u$ at $x_0$ equals $f$. Passing to the pointwise limit almost everywhere in the inequality from the Step 6 (possibly passing to a subsequence), we obtain that in a ball $B(x_0, \rho) \subset \Omega'$ we have
$$ f(x_0) - \delta \leq u(x) \leq f(x_0) + \delta.$$
As for arbitrary $\delta > 0$ there exists a ball $B(x_0, \rho)$ for which the above inequality is satisfied, we see that
$$ \lim_{\rho \rightarrow 0}\, \,  \text{ess sup}_{B(x_0,\rho)} |u(y) - f(x_0)| = 0,$$
so $Tu(x_0) = f(x_0)$. This equality holds for $\mathcal{H}^{N-1}$-almost all $x_0 \in \partial\Omega$, so $Tu = f$. \qed
\end{dd}

We conclude this Section with an existence result for norms such that $B_\phi(0,1)$ is not strictly convex. As we depend on an existence result in \cite{Gor3}, we will assume that $\Omega \subset \mathbb{R}^2$.

\begin{stw}\label{stw:existence2d}
Let $\phi$ be any norm on $\mathbb{R}^2$ and let $\Omega \subset \mathbb{R}^2$ be a strictly convex set with Lipschitz boundary. Suppose that $L^1(\partial\Omega)$ is a function such that $\mathcal{H}^{N-1}$-almost all points of $\partial\Omega$ are continuity points of $f$. Then there exists a minimiser to the Problem (\ref{aproblem}) with boundary data $f$.
\end{stw}

\begin{dd}
1. The case when $B_\phi(0,1)$ is strictly convex is already covered; it suffices to consider the case when $\partial B_\phi(0,1)$ has flat parts on the boundary.
Define $\phi_n = \phi + \frac{1}{n} l^2$. Then $B_{\phi_n}(0,1)$ is strictly convex and by Theorem \ref{thm:discexistencev1} there exists a solution $u_n \in BV(\Omega)$ to Problem (\ref{aproblem}) with boundary data $f$ with respect to the anisotropic norm $\phi_n$. In particular, $Tu_n = f$. Furthermore, the family $u_n$ is uniformly bounded in $BV(\Omega)$: we use Lemmata \ref{lem:poincare} and \ref{lem:variationestimate} respectively in the first and second inequality below to obtain
$$ \| u_n \|_{BV(\Omega)} = \int_\Omega |u_n|dx + \int_\Omega |Du_n| \leq C(\Omega) (\int_\Omega |Du_n| + \int_{\partial\Omega} |Tu_n| d\mathcal{H}^1) \leq $$
$$ \leq C(\Omega) (\int_{\partial\Omega} \phi_n(\nu) |Tu_n| d\mathcal{H}^1 + \int_{\partial\Omega} |Tu_n| d\mathcal{H}^1) \leq C(\Omega) \int_{\partial\Omega} (1 + \sup_{\partial B(0,1)} \phi + \frac{1}{n}) |f| d\mathcal{H}^{1} \leq M.$$
Hence, $u_n$ admits a subsequence (still denoted by $u_n$) convergent in $L^1(\Omega)$. By \cite[Theorem 3.2]{Gor3}, which is an anisotropic variant of Miranda's theorem when the anisotropic norm changes with $n$, $u$ is a function of $\phi-$least gradient; if we prove additionally that $Tu = f$, then $u$ is a solution to Problem (\ref{aproblem}).

2. Fix any continuity point $x_0 \in \partial\Omega$. Take arbitrary $\varepsilon > 0$. As $f \in C(\partial\Omega)$, there exists a neighbourhood of $x$ in $\partial\Omega$ such that 
$$ f(x) - \varepsilon \leq f(y) \leq f(x) + \varepsilon \qquad \text{ in } B(x, \delta_1) \cap \partial\Omega.$$
As $\Omega$ is strictly convex, for sufficiently small $\delta_1$ the set $B(x, \delta_1) \cap \partial\Omega$ consists of two points $p_1, p_2$ and the line segment $\overline{p_1 p_2}$ lies inside $\Omega$. Denote by $\Delta$ the open set bounded by an arc of $\partial\Omega$ containing $x$ and the line segment $\overline{p_1 p_2}$. Let us take a ball $B(x, \delta_2)$ such that $B(x, \delta_2) \cap \Omega \subset \Delta$. We argue as in Step 7 of the proof of Theorem \ref{thm:discexistencev1} that for every $n$ we have
$$ f(x) - \varepsilon \leq u_n(y) \leq f(x) + \varepsilon \qquad \text{ in } B(x, \delta_2) \cap \Omega;$$

3. As $u_n \rightarrow u$ in $L^1(\Omega)$, on some subsequence (still denoted $u_n$) we have convergence almost everywhere; hence
$$ f(x) - \varepsilon \leq u(y) \leq f(x) + \varepsilon \qquad \text{ for a.e. } y \in B(x, \delta_2) \cap \Omega.$$
As for arbitrary $\delta > 0$ there exists a ball $B(x_0, \rho)$ for which the above inequality is satisfied, we see that
$$ \lim_{\rho \rightarrow 0}\, \,  \text{ess sup}_{B(x_0,\rho)} |u(y) - f(x_0)| = 0,$$
so $Tu(x_0) = f(x_0)$. This equality holds for $\mathcal{H}^{N-1}$-almost all $x_0 \in \partial\Omega$, so $Tu = f$. \qed
\end{dd}

\subsection{Examples}

Now, we illustrate the results of this Section with a few examples. The first example gives yet another proof of an already known result, see \cite[Theorem 1.1]{Gor1}, \cite[Proposition 5.1]{DS} and \cite[Theorem 1.2]{RS}. This result serves as a toy problem and motivation for Theorem \ref{thm:discexistencev1}. Moreover, it is generalised to the anisotropic setting.

\begin{prz}
Let $\Omega \subset \mathbb{R}^2$ be a strictly convex set. Suppose that $f \in BV(\partial\Omega)$. Then the discontinuity set of $f$ is at most countable, hence it has $\mathcal{H}^1-$measure zero and by Theorem \ref{thm:discexistencev1} there exists a minimiser to Problem (\ref{problem}).  Furthermore, this reasoning extends to norms other than the Euclidean norm and to metric integrands satisfying (\ref{hypothesis}), provided that $\Omega$ satisfies the barrier condition.
\end{prz}

It is important to note that in dimension two the class of boundary data for which we have existence of minimisers is much larger than $BV(\partial\Omega)$; for instance, it contains all functions such that the set of discontinuity points is countable, even if the sum of jumps is infinite. However, there exist functions in $L^\infty(\partial\Omega)$ not covered by Theorem \ref{thm:discexistencev1}, for instance characteristic functions of fat Cantor sets; in particular, Theorem \ref{thm:discexistencev1} does not apply for the counterexample introduced in \cite{ST}. 

\begin{prz}
Let $\Omega = B(0,1) \subset \mathbb{R}^2$ and let $f = \chi_C$, where $C$ is the fat Cantor set constructed in \cite{ST}. This function is discontinuous on $C$, which has positive $\mathcal{H}^1-$measure, so Theorem \ref{thm:discexistencev1} does not guarantee existence of a minimiser to Problem (\ref{problem}); indeed, the authors of \cite{ST} prove that there is no minimiser.
\end{prz}

Another issue is the uniqueness of minimisers. A well-known example, attributed to John Brothers (see \cite[Example 2.7]{MRL}), shows that if the boundary data is discontinuous, then even in the isotropic case we may not expect uniqueness of minimisers.

\begin{prz}
Take boundary data given by the formula
$$ h(x,y) = \twopartdef{x^2 - y^2 + 1}{\text{ if } |x| > \frac{1}{\sqrt{2}}}{x^2 - y^2 - 1}{\text{ if } |x| < \frac{1}{\sqrt{2}}.}$$
Fix any $\lambda \in [-1,1]$. Then $u \in BV(\Omega)$ is a function of least gradient if and only if
$$ u(x,y) = \threepartdef{2x^2}{\text{ if } |x| > \frac{1}{\sqrt{2}}}{\lambda}{\text{ if } |x| < \frac{1}{\sqrt{2}}, |y| < \frac{1}{\sqrt{2}}}{-2y^2}{\text{ if } |y| > \frac{1}{\sqrt{2}}.}$$
The structure of level sets is precisely the same as for boundary data $h_0(x,y) = x^2 - y^2$, the only difference being that the values of the minimiser $u_0$ for boundary data $h_0$ are shifted in two opposite directions in different subsets of $\Omega$. It turns out that in the square where the function $u_0$ is constant we can choose freely the value of the shift. The situation is presented on Figure \ref{fig:brothers}. 

\begin{figure}[h]
    \includegraphics[scale=0.25]{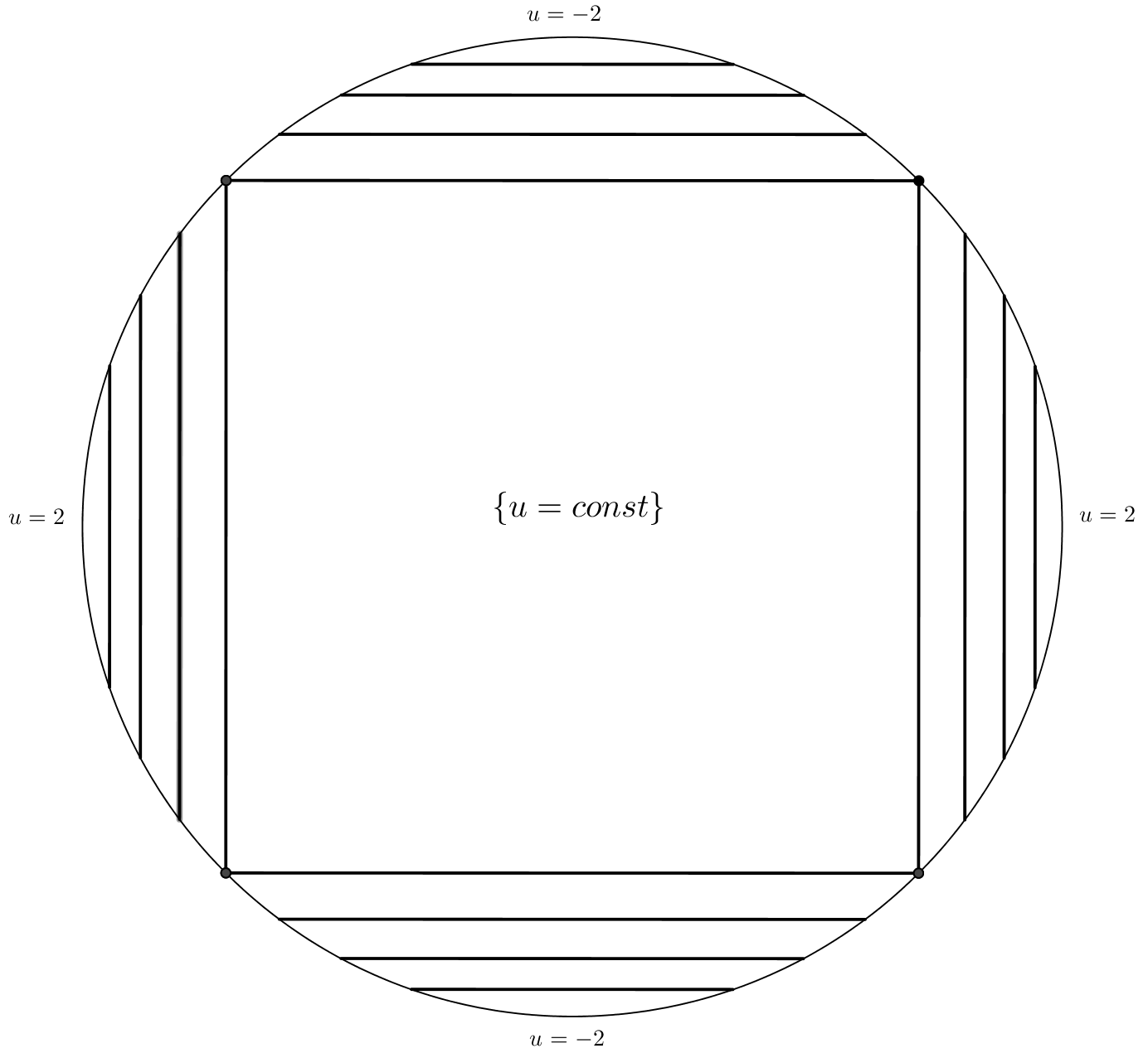}
    \caption{Nonuniqueness of minimisers for discontinuous boundary data}
    \label{fig:brothers}
\end{figure}

In general, the nonuniqueness is related to the formation of level sets of $u$ of positive Lebesgue measure, see \cite[Theorem 1.1]{Gor2}. The boundary data in this example has only four discontinuity points - it is the smallest number of discontinuity points for which we may lose uniqueness (the situation with fewer discontinuity points is considered for instance in \cite[Corollary 3.2]{GRS}).
\end{prz}

Finally, we abandon $\mathbb{R}^2$ and use Theorem \ref{thm:discexistencev1} to obtain existence of minimisers in a higher dimension. Moreover, the method used in Theorem \ref{thm:discexistencev1} is constructive and enables us to directly find the minimiser (provided that we can directly compute the result of Sternberg-Williams-Ziemer construction for the approximation).

\begin{prz}
Let $\Omega = B(0,1) \subset \mathbb{R}^3$. Take the boundary data to be
$$ f(x,y,z) = \twopartdef{1}{\text{ if } |z| > a}{-1}{\text{ if } |z| < a,}$$
where $a \in (0,1)$. The boundary data is continuous everywhere except the two circles which are intersections of $\partial\Omega$ and the planes $\{ z = \pm a \}$. Hence, the discontinuity set has $\mathcal{H}^2-$measure zero and there exists a minimiser to Problem (\ref{problem}).

Using an approximation as in the proof of Theorem \ref{thm:discexistencev1} and the axial symmetry, we can obtain the structure of minimisers, which differs with $a$. If $a$ is chosen so that the area of the two discs which are intersections of $\Omega$ and $\{ z = \pm a \}$ is smaller than the area of the part of the catenoid spanned by them, we have
$$ u(x,y,z) = \twopartdef{1}{\text{ if } |z| > a}{-1}{\text{ if } |z| < a.}$$
If the area of the discs is larger, then
$$ u(x,y,z) = \threepartdef{1}{\text{ if } |z| > a}{1}{\text{ if } |z| < a, \text{ inside the catenoid}}{-1}{\text{ if } |z| < a, \text{ outside the catenoid}.} $$
Finally, if the area of the discs is equal to the area of the part of catenoid spanned by them, mimimisers are no longer unique and
$$ u(x,y,z) = \threepartdef{1}{\text{ if } |z| > a}{\lambda}{\text{ if } |z| < a, \text{ inside the catenoid}}{-1}{\text{ if } |z| < a, \text{ outside the catenoid},}$$
where $\lambda \in [-1,1]$. By \cite[Theorem 1.1]{Gor2} these are all possible minimisers. Moreover, we observe that we can obtain all the minimisers in the last case using slightly different approximations, i.e. if we take the mollification kernel $\rho$ to be asymmetric.
\end{prz}

\section{Results for unbounded domains}\label{sec:unbounded}

In this Section, we consider the case when the domain $\Omega$ is unbounded. As in the least gradient problem for bounded domains, our main interest is to find what conditions do we need to impose on the domain and the boundary data in order to obtain existence and uniqueness of minimisers. Of course we will need to modify our notion of a solution; as we will see later, the solutions we construct need not lie in $BV(\Omega)$, but rather in $BV_{loc}(\Omega)$. Moreover, we have two kinds of additional difficulties: regularity of boundary data and shape of the domain.

For clarity, in this Section we will present the reasoning in the setting of the isotropic least gradient problem; we will remark when an analogous reasoning works also in the anisotropic case. Most notably the difference will appear with respect to the barrier condition and arguments concerning uniqueness. Following Miranda, see \cite{Mir}, we introduce the following definition of least gradient functions in an unbounded set $\Omega$:

\begin{dfn}
We say that $u \in BV_{loc}(\Omega)$ is a function of least gradient in $\Omega$ if for every function $g \in BV_{loc}(\Omega)$ with compact support $K \subset \Omega$ we have
$$ \int_{K} |Du| \leq \int_K |D(u+g)|.$$
We say that $u$ solves the least gradient problem on $\Omega$ with respect to boundary data $f \in L^1_{loc}(\partial\Omega)$, if both of the following conditions hold:
\begin{equation}\label{uproblem}\tag{ULGP}
u \text{ is a least gradient function in } \Omega \quad and 
\end{equation}
\begin{equation*}
\text{for } \mathcal{H}^{N-1}\text{-a.e. } x \in \partial\Omega \text{ we have } \dashint_{B(x,r) \cap \Omega} |u(y) - f(x)| \rightarrow 0 \text{ as } r \rightarrow 0.
\end{equation*}
For bounded domains this is equivalent to Problem (\ref{problem}); however, here we cannot minimise the total variation, as the total variation turns out to be finite only in the most restrictive cases. Furthermore, the trace condition is understood pointwise, $\mathcal{H}^{N-1}$-a.e. on the boundary of $\Omega$. As the boundary of $\Omega$ is unbounded, continuous functions need not be bounded; this is our first additional difficulty:
\end{dfn}

{\it A. Regularity of boundary data.} In this paper, we are mostly interested in the case of continuous boundary data. However, for unbounded domains the Sternberg-Williams-Ziemer construction (see \cite{SWZ}) does not work - the construction relies on two auxiliary problems, involving minimalisation of perimeter and maximisation of area, which do not need to have solutions if the domain is unbounded. For this reason, we are going to consider multiple function spaces on $\partial\Omega$ and make distinctions between them when it comes to existence, uniqueness and regularity of minimisers. We consider the following spaces:

1. $L^1(\partial\Omega) \cap C_0(\partial\Omega)$. This is the most natural function space to arise in this problem - for instance, we may regard the trace as an operator on $\Omega$ and not merely as a pointwise property that holds almost everywhere.

2. $C_0(\partial\Omega)$. This class arises naturally if we try to approximate the boundary data with continuous boundary data with compact support, which is the approximation that yields uniqueness of minimisers.

3. $C_b(\partial\Omega)$ and $C(\partial\Omega)$. In this class there appear some interesting phenomena concerning the shape of superlevel sets, such as creation of ''shock waves'' that extend to infinity. In particular, this leads to nonuniqueness of minimisers for a wide class of domains.

4. The case when the data are continuous almost everywhere. In particular, this covers the case $BV(\partial\Omega)$ for $\Omega \subset \mathbb{R}^2$. As we know from the bounded domain case, we cannot expect much more than existence of minimisers - this case combined with the results from Section \ref{sec:discontinuous} is a corollary to the previous ones.

We cannot expect much less regularity; in the case when $\Omega$ is a disk, see \cite{ST} for an example of a function that is $L^\infty(\Omega)$, is discontinuous on a set of positive measure, and there is no minimiser. The example adapts well to the unbounded case. \\

{\it B. Shape of domains.} The second kind of assumptions concerns the shape of domains. As in the bounded case, we have to assume a condition similar to strict convexity of the domain. We will see that the shape of the domain makes no difference in the existence proof; however, the shape of the domain may influence the regularity of the resulting minimiser. We are interested in three kinds of domains: 

1. Strictly convex unbounded domains such that $\Omega \neq \mathbb{R}^N$. In particular, $\Omega$ lies in a halfspace. In this class, we are able to obtain existence of minimisers and uniqueness of minimisers for data in $C_0(\partial\Omega)$.

2. Domains with special features: firstly, domains which are in a sense ''one-dimensional'', i.e. domains that are unbounded only in one direction. In particular, these domains lie in a strip, so we have a uniform Poincar\'e inequality for any open subset $\Omega' \subset \Omega$ with Lipschitz boundary. Secondly, we will consider domains which contain a cone; these domains will be crucial to the phenomenon of nonuniqueness.

3. Finally, we want to consider boundary values in infinity. This is motivated mostly by the case when the domain $\Omega$ contains a cone. For simplicity, we restrict ourselves to the whole $\mathbb{R}^N$ and consider a standard compactification of $\mathbb{R}^N$ defined by adding a point in each direction; the resulting space is denoted by $\mathbb{R}^N \cup \partial B(0,1)$ and the boundary data are given as a function $f \in L^1(\partial B(0,1))$.

\subsection{Existence of minimisers}

The Theorem below proves existence of minimisers for the least gradient problem on an unbounded domain in full generality. Later, we will consider what modifications can we make to the proof below in order to obtain some additional regularity or uniqueness of minimisers.

\begin{tw}\label{thm:generalexistence}
Let $\Omega \subset \mathbb{R}^N$ be a strictly convex set and $\Omega \neq \mathbb{R}^N$. Let $f \in C(\partial\Omega)$. Then there exists a minimiser $u \in BV_{loc}(\Omega)$ of the least gradient problem with boundary data $f$.
\end{tw}

\begin{dd}
1. We begin by noting that as $\Omega$ is a convex set which is not equal to $\mathbb{R}^N$, it is contained in a halfspace; it suffices to fix any $x \in \partial\Omega$ and consider any supporting hyperplane $H$ with inward normal vector $\nu$. Then $\Omega$ lies entirely on one side of this hyperplane, which is a halfspace we denote by $H_+$. The other halfspace $H_-$ is disjoint with $\Omega$. We notice that the shifted halfspaces of the form $H_- + t\nu$ for $t > 0$ intersect $\Omega$ and their union is $\mathbb{R}^N$, so they cover the whole $\Omega$.

2. Now, we introduce both approximating sets $\Omega_n$ and the approximate boundary data $f_n \in C(\partial\Omega_n)$. Choose $M_n$ to be a sequence of positive numbers such that $M_n \geq M_{n-1} + 2$. We set $f_n$ to be a continuous function with compact support such that $f_n = f$ in $\partial\Omega \cap (H_- + M_n \nu)$, $f_n = 0$ in $\partial\Omega \backslash (H_- + (M_n + 1)\nu)$ and (by a variant of Tietze extension theorem) as a continuous function with values in the line segment $[- f(x), f(x)]$ in $\partial\Omega \cap ((H_- + (M_n + 1)\nu) \backslash (H_- + M_n \nu))$. In particular, the sequence $f_n$ converges to $f$ locally uniformly and the sequence $|f_n|$ is monotone and converges locally uniformly to $|f|$.

3. Let $\Omega_n$ be an increasing sequence of strictly convex sets such that $\bigcup_{m = 1}^\infty \Omega_n = \Omega$ and such that $\Omega_n \cap (H_- + M_n \nu) = \Omega \cap (H_- + M_n \nu)$. The construction is shown on Figure \ref{fig:konstrukcja}; the shaded set is $\Omega_1$. We may consider $f_n$ to be a function on $\Omega_n$ via a simple identification: let $\widetilde{f}_n \in C(\partial\Omega_n)$ be defined as $\widetilde{f}_n = f_n$ on $\partial\Omega_n \backslash \partial\Omega$ and $\widetilde{f}_n = 0$ on $\partial\Omega_n \backslash \partial\Omega$. Using the Sternberg-Williams-Ziemer construction we obtain a minimiser $u_n \in BV(\Omega_n)$ of the least gradient problem in $\Omega_n$ with boundary data $\widetilde{f}_n$. By \cite[Proposition 4.1]{GRS} the restriction of $u_n$ to $\Omega_m$ with $m \leq n$ also lies in $BV(\Omega_m)$ and is a function of least gradient.

\begin{figure}[h]
    \includegraphics[height=5cm, width=7cm]{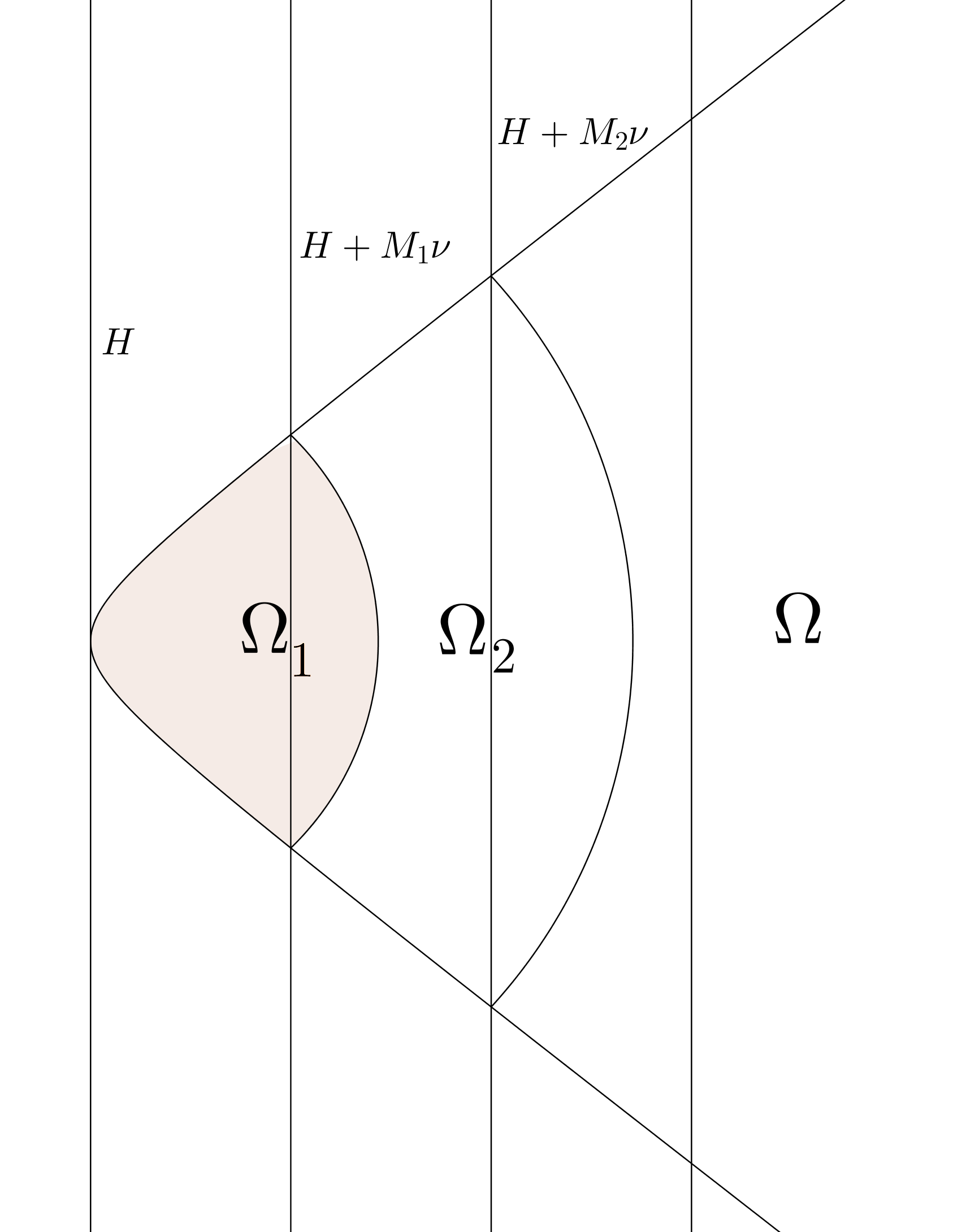}
    \caption{The construction of the approximating sets $\Omega_n$}
    \label{fig:konstrukcja}
\end{figure}

4. For every $m$, we need to show a uniform bound in $BV(\Omega_m)$ for the sequence of approximations in order to obtain existence of a limit function in the topology of $L^1(\partial\Omega)$. Let $n > m$ and consider the restrictions $u_n |_{\Omega_m}$. We calculate
$$ \| u_n \|_{BV(\Omega_m)} = \int_{\Omega_m} |u_n| + \int_{\Omega_m} |Du_n| \leq C(\Omega_m) (\int_{\Omega_m} |Du_n| + \int_{\partial \Omega_m} |Tu_n|) \leq$$
$$ \leq 2 C(\Omega_m) \int_{\partial \Omega_m} |Tu_n| = 2 C(\Omega_m) (\int_{\partial \Omega_m \cap \partial\Omega} |Tu_n| + \int_{\partial \Omega_m \backslash\partial\Omega} |Tu_n|) \leq $$
$$ = 2 C(\Omega_m) (\int_{\partial \Omega_m \cap \partial\Omega} |f| + \int_{\partial \Omega_m \backslash\partial\Omega} \sup_{\Omega_m} |u_n|) \leq 2\widetilde{C}(\Omega_m) \sup_{\partial\Omega \cap (H_- + M_{n+1} \nu)} |f|.$$
Hence the sequence $u_n|_{\Omega_m}$ is bounded in $BV(\Omega_m)$ and admits a convergent subsequence in $L^1(\Omega_m)$ and almost everywhere. Using a diagonal argument, we obtain existence of of a subsequence $u_{n_k}$ (extended by zero outside $\Omega_{n_k}$) convergent to $u \in BV_{loc}(\Omega)$ in $L_{loc}^1(\Omega)$ and almost everywhere. By Theorem \ref{thm:miranda} $u$ is a function of least gradient in $\Omega$.

5. We have to ensure that $Tu = f$. We proceed as in \cite{Gor3}; as $f$ is continuous, we fix any point $x_0 \in \partial\Omega$ such that the mean integral condition in the definition of trace is satisfied. For every $\delta > 0$ we find a ball $B(x_0, r(\delta))$ such that every $u_n$ satisfies 
$$ f(x_0) - \delta \leq u_n(x) \leq f(x_0) + \delta \qquad \text{ in } B(x_0, r(\delta)) \cap \Omega.$$
Hence, as in the proof of Theorem \ref{thm:discexistencev1}, $u$ satisfies the same inequalities in $\Omega \cap B(x_0, r'(\delta))$ and so $Tu(x_0) = f(x_0)$. As this holds for almost all $x_0 \in \partial\Omega$, $Tu = f$. \qed
\end{dd}

\begin{wn}
Let $\Omega \subset \mathbb{R}^N$ be a strictly convex set and $\Omega \neq \mathbb{R}^N$. Let $f \in L^\infty_{loc}(\partial\Omega)$ such that $\mathcal{H}^{N-1}$-almost all points of $\partial\Omega$ are continuity points of $f$. Then there exists a minimiser $u \in BV_{loc}(\Omega)$ of the least gradient problem with boundary data $f$.
\end{wn}

\begin{dd}
The proof of Theorem \ref{thm:generalexistence} requires only a few modifications. We choose the sequence $\Omega_n$ in the same way and construct the approximating sequence $f_n$ in a simpler way: we extend $f$ by zero on $\partial \Omega_m \backslash \partial\Omega$. The resulting function $f_n \in L^\infty(\partial\Omega_n)$ is such that $\mathcal{H}^{N-1}$-almost all points of $\partial\Omega$ are continuity points of $f$. By Theorem \ref{thm:discexistencev1} there exists a minimiser to Problem (\ref{problem}) $u_n \in BV(\Omega_n)$. Now, we notice that the uniform estimates in Step 4 depend only on the local $L^\infty$ bounds on $f|$ and not on its continuity; finally, we repeat Step 5 only for continuity points of $f$. \qed
\end{dd}

\begin{uw}
The results above hold with minor changes also for norms other than the Euclidean norm; the greatest difficulty is to ensure that the barrier condition is satisfied. For instance, we could use the translation invariance and obtain $\Omega_m$ by cutting $\Omega$ using a hyperplane as in the proof of Theorem \ref{thm:generalexistence} and then reflecting the resulting set with respect to this hyperplane. However, when $\phi$ depends on location, this is not immediate and depends on the exact form of $\phi$. If $\partial\Omega \in C^2$, this boils down to solving a degenerate elliptic equation, see \cite[Lemma 3.2]{JMN}. Moreover, we need condition (\ref{hypothesis}) to hold in order to use a comparison principle as in the proof of Theorem \ref{thm:discexistencev2} to prove that the trace of the minimiser is correct.
\end{uw}

\subsection{Regularity}

In this subsection, we discuss the main regularity features of minimisers in the least gradient problem on unbounded domains. The main result in this subsection, Proposition \ref{prop:continuity}, concerns continuity of minimisers. While minimisers to Problem (\ref{problem}) for continuous boundary data are continuous up to the boundary if the domain $\Omega$ is bounded, it is not necessarily obvious here - our approximation procedure provides gives no estimates which enable us to prove that the convergence $u_n \rightarrow u$ is locally uniform. Furthermore, this approach would only prove continuity of the one minimiser obtained using the approximation procedure; as we will see in Example \ref{ex:nonuniqueness}. Proposition \ref{prop:continuity} goes around both of these problems. In order to prove it, we first recall a version of maximum principle for minimal surfaces.

\begin{stw}\label{swz:zasadamaksimum}
(\cite[Theorem 2.2]{SWZ}) Suppose that $E_1 \subset E_2$ and let $\partial E_1, \partial E_2$ be area-minimising in an open set $U$. Further, suppose that $x \in \partial E_1 \cap \partial E_2 \cap U$. Then $\partial E_1$ and $\partial E_2$ coincide in a neighbourhood of $x$.
\end{stw}

\begin{stw}\label{prop:continuity}
Let $\Omega \subset \mathbb{R}^N$ be a strictly convex set and $\Omega \neq \mathbb{R}^N$. Let $f \in C(\partial\Omega)$. Suppose that $u \in BV_{loc}(\Omega)$ is a minimiser of the least gradient problem with boundary data $f$. Then $u \in C(\overline{\Omega})$.
\end{stw}

For a function $u \in BV_{loc}(\Omega)$ we say that $x \in S_u$, the approximate discontinuity set of $u$, if the lower and upper approximate limits $u^{\wedge}(x)$ and $u^{\vee}(x)$ do not coincide (see for instance \cite[Chapter 3]{AFP}).

\begin{dd}
Let $x \in \Omega$ be a point such that $u$ is not continuous at $x$. By \cite[Theorem 4.1]{HKLS} we have $x \in S_u$. By the definition of $S_u$, $x \in \partial \{ u \geq t \}$ for each $t \in (u^{\vee}(x), u^{\wedge}(x))$. As for $s > t$ we have $\{ u \geq s \} \subset \{ u \geq t \}$, by Proposition \ref{swz:zasadamaksimum} the connected components of $\partial \{ u \geq t  \}$ for $t \in (u^{\vee}(x), u^{\wedge}(x))$ passing through $x$ agree; we will denote this surface by $S$.

We will see that $S \cap \partial \Omega = \emptyset$. Suppose the contrary and take $z_0 \in S \cap \partial \Omega$. Fix $\varepsilon > 0$. Then, by continuity of $f$, in a neighbourhood of $z_0$ we have
$$ f(z_0) - \varepsilon \leq f(z) \leq f(z_0) + \varepsilon.$$
Using the same argument as in Steps 6 and 7 of the proof of Theorem \ref{thm:discexistencev1}, we obtain that for $y \in B(z_0, r(\varepsilon)) \cap \Omega$ we have
$$ f(z_0) - \varepsilon \leq u(y) \leq f(z_0) + \varepsilon.$$
However, as $z_0 \in S \cap \partial\Omega$, in any neighbourhood of $z_0$ in $\Omega$ there are values of $u$ greater or equal to $u^{\wedge}(x)$ and smaller or equal to $u^{\vee}(x)$. This contradicts the estimates above for sufficiently small $\varepsilon$.

As $S \cap \partial \Omega = \emptyset$ and $S$ is a minimal surface, $S$ extends to infinity; but then we could replace $S$ by a truncated $S$ (in the notation of the proof of Theorem \ref{thm:generalexistence}, using a projection of $S$ onto a plane $H + s\nu$) and obtain a surface with lower area. Hence for $t \in (u^{\vee}(x), u^{\wedge}(x))$ the set $\{ u\geq t \}$ is not minimal, which contradicts Theorem \ref{tw:bgg}.

Now, take $x \in \partial\Omega$. Then, as $u$ is continuous inside $\Omega$, the essential supremum is the same as supremum and we see that
$$ \lim_{r \rightarrow 0} \text{ess sup}_{B(x,r)} |u(y) - f(x)| = \lim_{r \rightarrow 0} \sup_{B(x,r)} |u(y) - f(x)| = 0,$$
hence $u$ is continuous at $x$. \qed
\end{dd}

A careful inspection of the proof above yields a generalisation of Lemma \ref{lem:swzbdry}:

\begin{wn}
Let $\Omega \subset \mathbb{R}^N$ be a strictly convex set and $\Omega \neq \mathbb{R}^N$. Suppose that $u \in BV_{loc}(\Omega)$ is a least gradient function with trace $f \in L^1_{loc}(\partial\Omega)$. Then for each $t \in \mathbb{R}$ we have an inclusion
$$ \partial \{ u \geq t \} \cap \partial \Omega \subset f^{-1}(t) \cup D,$$
where $D$ is the set of discontinuity points of $f$.
\qed
\end{wn}

\begin{uw}
In the anisotropic setting, the assumptions concerning $\phi$ in the continuity proof are more restrictive than in the existence proof. The main problem is that we use Proposition \ref{tw:bgg}. Hence, Proposition \ref{prop:continuity} is valid in the anisotropic case in two settings: firstly, when $\Omega \subset \mathbb{R}^2$ and $\phi$ is a norm such that $B_\phi(0,1)$ is strictly convex, the connected components of $\phi-$area-minimising sets are line segments and we do not need to use the maximum principle. Secondly, when an anisotropic version of the maximum principle holds, for instance for weighted least gradient problem, when $\phi(x, Du) = a(x)|Du|$ with $a \in C^2(\overline{\Omega})$, see \cite[Theorem 3.1]{Zun}.
\end{uw}

However, contrary to the results in the bounded domain case, we cannot expect H\"older continuity of minimisers. This is due to the fact that $\partial \Omega$ becomes asymptotically flat as $|x| \rightarrow \infty$; it means that the regularity of minimisers at infinity is the same as near a point where mean curvature vanishes and its growth rate is slower than polynomial. For the fact that in the neighbourhood of such points we may lose H\"older continuity, see \cite[Remark 5.8]{SWZ}.

Now, we will see that integrability of the minimiser and its total variation can only happen under very special circumstances, both in terms of the regularity of boundary data and the shape of the domain. 

\begin{stw}
Let $\Omega \subset \mathbb{R}^N$ be a domain that is unbounded only in one direction and such that its cross-sectional area is uniformly bounded. Suppose that $u \in BV_{loc}(\Omega)$ is a minimiser to the least gradient problem with boundary data $f \in L^1(\partial\Omega) \cap L^\infty(\Omega)$. Then $u \in BV(\Omega)$.
\end{stw}

\begin{dd}
We are going to utilise the Poincar\'{e} inequality. Firstly, let us see that the Poincar\'{e} inequality holds in $\Omega$ for each $\Omega_m$, which is $\Omega$ cut at level $m$ by a hyperplane, uniformly in $m$ - the constant in the inequality depends only on the width $d$ of the strip. 

Now, we estimate
$$ \int_{\Omega_m} |u| \leq C(d) (\int_{\Omega_m} |Du| + \int_{\partial \Omega_m} |Tu|)$$
and so
$$ \| u \|_{BV(\Omega_m)} \leq C(d) \int_{\partial \Omega_m} |Tu| + (C(d) + 1) \int_{\Omega_m} |Du| \leq (2C(d) + 1) \int_{\partial \Omega_m} |Tu| \leq  $$
$$ \leq (2C(d) + 1) (\int_{\partial\Omega} |f| + \int_{\Gamma_m} |Tu|) \leq (2C(d) + 1) (\| f \|_{L^1(\partial\Omega)} + d^{N-1} \sup_{\partial\Omega} |f|). $$
Hence the $BV$ norm of $u$ is uniformly bounded in each $\Omega_m$ by quantities which only depend on the width $d$ (in all but one dimensions), the supremum of $|f|$ and the $L^1$ norm of $f$. Thus
$$ \| u \|_{BV(\Omega)} \leq \widetilde{C}(d) (\| f \|_{L^1(\partial\Omega)} + \sup_{\partial\Omega} |f|).$$
We point out that this proof only did not require continuity of $f$, only the boundedness. \qed
\end{dd}

If $f \notin L^1(\partial\Omega)$ or $\Omega$ is not one-dimensional, then we cannot hope that the minimiser $u$ (if it exists) is in $BV(\Omega)$; let us see two simple examples:

\begin{prz}
Let $\Omega \subset \mathbb{R}^2$ be defined as
$$ \Omega = \{ (x,y): \, x > 0, \, e^{-x} - 1 \leq y \leq -e^{-x} + 1 \}.$$
Now, we will construct the boundary data $f \in C_0(\partial\Omega)$. Let $\rho$ be a standard mollifier on $\mathbb{R}$ with support in $(-1,1)$. Now, let $f$ be defined by the formula

$$ f(x,y) = \sum_{n=1}^\infty \frac{1}{n} \rho(x - n^2).$$
Then the minimiser to the least gradient problem given by Theorem \ref{thm:generalexistence} is
$$ u(x,y) = \sum_{n=1}^\infty \frac{1}{n} \rho(x - n^2),$$
i.e. all level lines are vertical. However, $u \notin L^1(\Omega)$ and the total variation of $u$ is infinite.
\end{prz}

\begin{prz}
Let $\Omega \subset \mathbb{R}^2$ be defined as
$$ \Omega = \{ (x,y): \, x > 0, \, |y| \leq x^3 \}.$$
Let $f(x,y) = \frac{1}{(x+1)^2}$. Then the minimiser in the least gradient problem exists and again all level lines are vertical and equals $u(x,y) = \frac{1}{(x+1)^2}$; however, $u \notin L^1(\Omega)$ and the total variation of $u$ is infinite.
\end{prz}

Finally, we will briefly discuss a new phenomenon that happens when $f \notin C_0(\partial\Omega)$ and is the reason behind the lack of uniqueness: there may form level sets which do not connect points from $\partial \Omega$, but instead escape to infinity. In particular, if $u$ is a solution of Problem (\ref{uproblem}), then there may exist connected components of $\partial \{ u \geq t \}$ with infinite area. However, we will see that in dimension two there is at most one such connected component.

\begin{figure}[h]
    \includegraphics[scale=0.4]{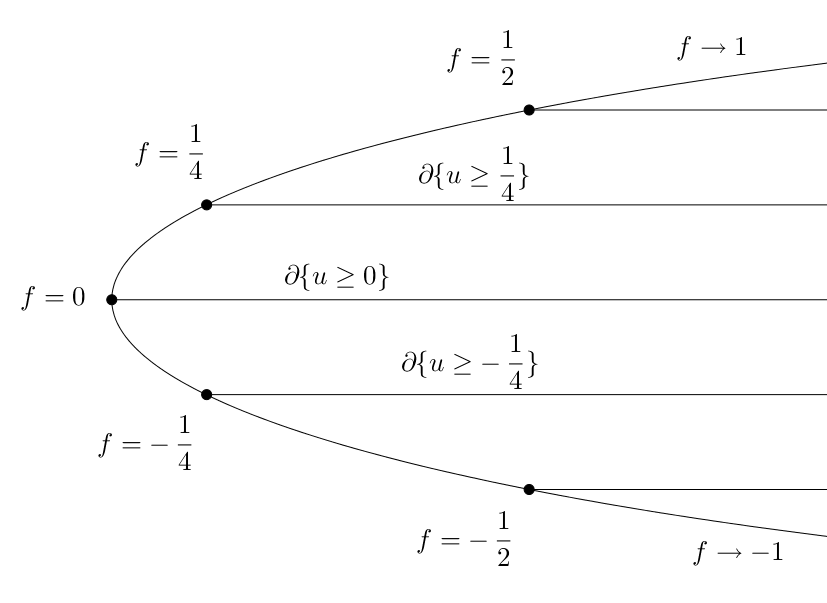}
    \caption{The level lines escape to infinity}
    \label{fig:monotone}
\end{figure}

\begin{prz}
Let $\Omega \subset \mathbb{R}^2$ be defined as
$$ \Omega = \{ (x,y): \, x > 0, \, e^{-x} - 1 \leq y \leq -e^{-x} + 1 \}.$$
We take the boundary data $f(x,y) = y$. Then on the lower branch of $\partial\Omega$ we have $f(x,y) \rightarrow -1$ as $(x,y)$ tends to infinity and on the upper branch of $\partial\Omega$ we have $f(x,y) \rightarrow 1$ as $(x,y)$ tends to infinity. Therefore for all $t \in (-1,1)$ the set $\partial\{ u \geq t \}$ is a halfline, going in the horizontal direction to the right, starting at a point of the form $(x,t)$. The situation is presented on Figure \ref{fig:monotone}.
\end{prz}

\begin{stw}
Let $\Omega \subset \mathbb{R}^2$ be an open unbounded strictly convex set. Let $u \in BV_{loc}(\Omega)$ be a least gradient function. Then, for each $t \in \mathbb{R}$ there is at most one unbounded connected component of $\partial \{ u \geq t \}$.
\end{stw}

\begin{figure}[h]
    \includegraphics[scale=0.3]{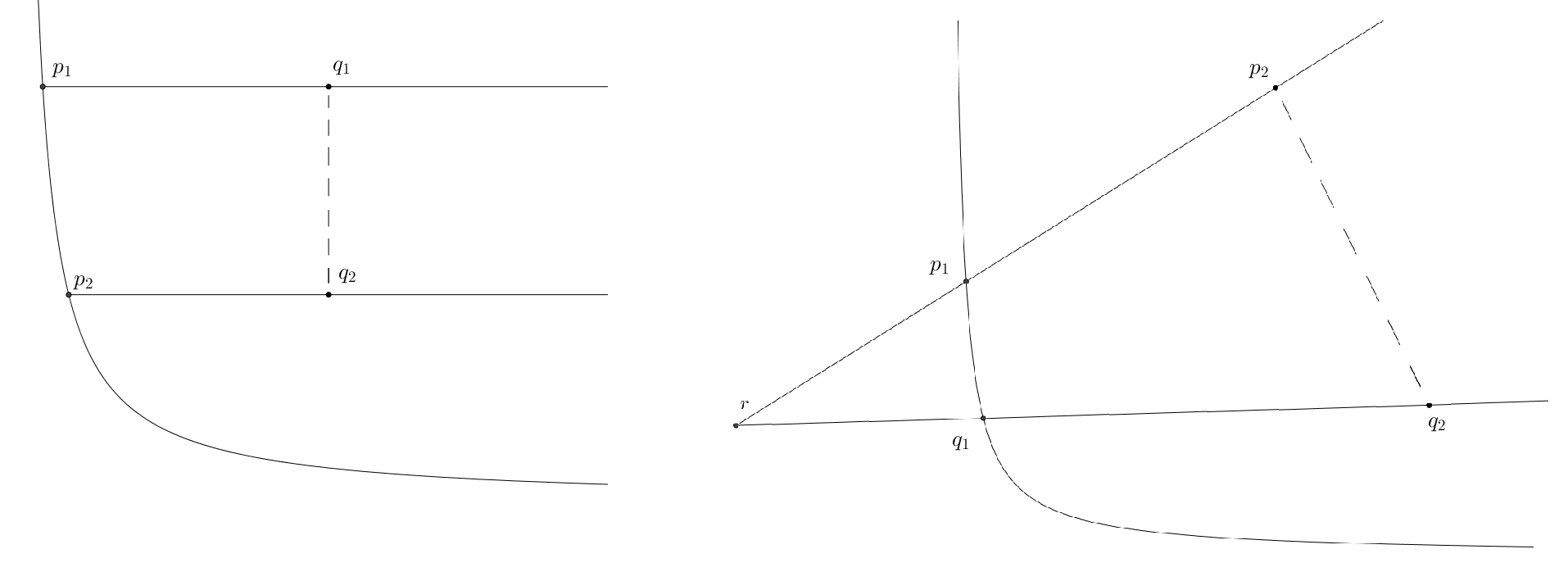}
    \caption{There is at most one level line escaping to infinity}
    \label{fig:escape}
\end{figure}

\begin{dd}
Suppose that $\partial \{ u \geq t \}$ has at least two unbounded connected components; in dimension two these are halflines with starting point on $\partial\Omega$. If these halflines are parallel, then we can replace them with a U-shaped polygonal chain consisting of two halflines and a line segment, locally reducing the total variation; the situation is presented on the left hand side of Figure \ref{fig:escape}. If we choose points $p_2, q_2$ sufficiently far from $\partial\Omega$, then the line segment $\overline{p_2 q_2}$ is shorter than the line segments $\overline{p_1 p_2}$ and $\overline{q_1 q_2}$ and hence the function $\chi_{\{ u \geq t \}}$ was not a function of least gradient, which contradicts Theorem \ref{tw:bgg}.

If these halflines are not parallel, they intersect in a point $r \notin \Omega$. Then we can replace them by another U-shaped polygonal chain if the line segment is far enough from $r$; the situation is presented on the right hand side of Figure \ref{fig:escape}. By the triangle inequality, the line segment $\overline{p_2 q_2}$ is shorter than the union of the line segments $\overline{p_2 r}$ and $\overline{r q_2}$. If we choose points $p_2, q_2$ sufficiently far from $\partial\Omega$, then the line segment $\overline{p_2 q_2}$ is also shorter than the union of the line segments $\overline{p_1 p_2}$ and $\overline{q_1 q_2}$. Hence the function $\chi_{\{ u \geq t \}}$ was not a function of least gradient, which contradicts Theorem \ref{tw:bgg}. \qed
\end{dd}

\subsection{Uniqueness of minimisers}

It turns out that the natural space for uniqueness of minimisers is the space $C_0(\partial\Omega)$. In this space, we can infer the uniqueness of minimisers from the uniqueness of minimisers in bounded domains in a similar way as in the existence proof, only using a more careful approximation. If the boundary data is less regular than $C_0(\partial\Omega)$, we may construct nonunique minimisers even for very simple boundary data.

\begin{tw}\label{thm:uniqueness}
Let $\Omega \subset \mathbb{R}^N$ be an unbounded strictly convex set and $\Omega \neq \mathbb{R}^N$. Let $f \in C_0(\partial\Omega)$. Then there exists a unique minimiser $u \in BV_{loc}(\Omega) \cap C(\overline{\Omega})$ of the least gradient problem with boundary data $f$.
\end{tw}

\begin{dd}
1. We recall a bit of the Sternberg-Williams-Ziemer construction of minimisers for continuous boundary data and a bounded strictly convex set $\Omega'$ (see \cite{SWZ}). Let $g \in C(\partial \Omega')$. We take its extension $G \in C(\mathbb{R}^N \backslash \Omega') \cap BV(\mathbb{R}^N \backslash \Omega')$ with compact support and denote $L_t = \{ G \geq t \}$. For almost all $t$, the superlevel sets $E_t = \{ u \geq t \}$ are solutions of the problem
$$ \min \{ P(E, \mathbb{R}^N): E \backslash \overline{\Omega'} = L_t \backslash \overline{\Omega'} \} $$
$$ \max \{ |E|: \, \, E \text{ is a solution of the above}\}.$$
The result does not depend on the choice of the extension $G$. In particular, the sets $E_t$ are defined uniquely.

2. As in the proof of Theorem \ref{thm:generalexistence}, fix any $x \in \partial\Omega$ and consider any supporting hyperplane $H$ with inward normal vector $\nu$. Let us take the halfspace $H_-$, which is disjoint with $\Omega$ and whose boundary is $H$. Again, the shifted halfspaces of the form $H_- + s\nu$ for $s > 0$ intersect $\Omega$ and their union is $\mathbb{R}^N$, so they cover the whole $\Omega$. As $f \in C_0(\partial\Omega)$, for every $n$ there exists $M_n$ so that 
$$ |f(x)| \leq \frac{1}{n} \qquad \text{ in } \partial\Omega \backslash (H_- + M_n \nu).$$
We may additionally require that $M_n \geq M_{n-1} + 2$. Let $\Omega_n$ be an increasing sequence of strictly convex sets such that $\bigcup_{m = 1}^\infty \Omega_n = \Omega$ and such that $\Omega_n \cap (H_- + M_n \nu) = \Omega \cap (H_- + M_n \nu)$. By the definition of $M_n$, for $t > \frac{1}{n}$ we have an inclusion 
$$\{ f \geq t \} \subset \partial\Omega \cap (H_- + M_n \nu).$$

3. By Proposition \ref{prop:continuity} the minimisers to the least gradient problem are continuous up to the boundary. Suppose that $u \in BV_{loc}(\Omega) \cap C(\overline{\Omega})$ and $v \in BV_{loc}(\Omega) \cap C(\overline{\Omega})$ are two minimisers to the least gradient problem on $\Omega$. Using an argument with projections as in Step 7 of the proof of Theorem \ref{thm:discexistencev1}, we obtain that
$$\{ u \geq t \} \subset \Omega \cap (H_- + M_n \nu) \qquad \text{and} \qquad \{ v \geq t \} \subset \Omega \cap (H_- + M_n \nu).$$
Let us consider the restrictions of $u$ and $v$ to $\overline{\Omega_n}$. Both functions are continuous, hence $u$ solves the least gradient problem on a bounded strictly convex domain $\Omega$ with boundary data
$$ \widetilde{u} = \twopartdef{f}{\text{ on } \partial\Omega \cap \partial \Omega_n}{u}{\text{ on } \partial \Omega_n \backslash \partial \Omega}$$
and analogously $v$ solves the least gradient problem on $\Omega_n$ for boundary data
$$ \widetilde{v} = \twopartdef{f}{\text{ on } \partial\Omega \cap \partial \Omega_n}{v}{\text{ on } \partial \Omega_n \backslash \partial \Omega.}$$

4. As $\widetilde{u}$ and $\widetilde{v}$ agree on $\partial\Omega \cap \partial \Omega_n$, we can choose their extensions $\widetilde{U}, \widetilde{V} \in C(\mathbb{R}^N \backslash \Omega_n) \cap BV(\mathbb{R}^N \backslash \Omega_n)$ to agree in a neighbourhood of $\partial\Omega \cap (H_- + M_n \nu)$. Moreover, by continuity the functions $\widetilde{U}, \widetilde{V}$ are less or equal to $t > \frac{1}{n}$ in a neighbourhood of $\partial\Omega \cap \partial\Omega_n$.

5. Pick $t$ such that both $\{ u \geq t \}$ and $\{ v \geq t \}$ solve the problem from the Sternberg-Williams-Ziemer construction for boundary data $\widetilde{u}$ and $\widetilde{v}$ respectively on $\partial \Omega_n$. The set $\{ u \geq t \}$ is determined only by the set $\{ \widetilde{U} \geq t \}$. However, it agrees with the set $\{ \widetilde{V} \geq t \}$ in a neighbourhood of $\Omega_n$, as $\widetilde{U} = \widetilde{V}$ in a neighbourhood of $\partial\Omega \cap (H_- + M_n \nu)$ and both sets are empty in a neighbourhood of $\partial\Omega \backslash (H_- + M_n \nu)$. By the uniqueness of the set resulting from the Sternberg-Williams-Ziemer procedure we have that $\{ u \geq t \} = \{ v \geq t \}$. Hence, for almost all $t > \frac{1}{n}$ we have $\{ u \geq t \} = \{ v \geq t \}$. We proceed similarly for $t < - \frac{1}{n}$. Hence almost all superlevel sets of $u$ and $v$ are uniquely determined and equal, so $u = v$ almost everywhere. \qed
\end{dd}

\begin{uw}
Clearly, the Proposition above holds also in the case when $f \in C(\partial\Omega)$ and $f$ has a finite limit $f_0$ as $|x| \rightarrow \infty$. Furthermore, it holds also if the limit is infinite; if $f \rightarrow + \infty$ as $|x| \rightarrow \infty$, then we have to choose the sequence $M_n$ in Step 1 so that
$$ f(x) \geq n \qquad \text{ in } \partial\Omega \backslash (H_- + M_n \nu)$$
and continue the proof as above.
\end{uw}

\begin{prz}\label{ex:nonuniqueness}
Let $\Omega \subset \mathbb{R}^2$ be defined as
$$ \Omega = \{ (x,y): \, x > 0, \, y > 0, \, \, xy > 1 \}. $$
Let the boundary data equal $f(x,y) = e^{-x}$. The boundary data is monotone, $f \in C_b(\partial\Omega)$ and it has finite limits in infinity.

Then the functions $u_1, u_2: \Omega \rightarrow \mathbb{R}$ defined by $u_1(x,y) = e^{-x}$ and $u_2(x,y) = e^{-\frac{1}{y}}$ with boundary data $f$. The first one has all level lines vertical and the second one has all level lines horizontal. Moreover, each value is attained only on one half-line; if we take any strictly convex $\Omega' \subset \Omega$, then $u_i$ restricted to $\partial \Omega'$ has only one minimum and one maximum. Hence both are functions of least gradient and the minimisers to Problem (\ref{uproblem}) are not unique. This situation is presented on Figure \ref{fig:nonuniqueness}. We note that we could obtain an uncountable family of minimisers by choosing the angles of incidence of the level lines. In particular, even though $f$ is monotone, we can obtain a least gradient function level set of positive (infinite) measure if we set for $x_0 \in \mathbb{R}_+$
$$ u_3 = \threepartdef{e^{-x}}{\text{ if } |x| < x_0}{e^{-\frac{1}{y}}}{\text{ if } |y| < \frac{1}{x_0}}{e^{-x_0}}{\text{ otherwise.}}$$

\begin{figure}[h]
    \includegraphics[scale=0.32]{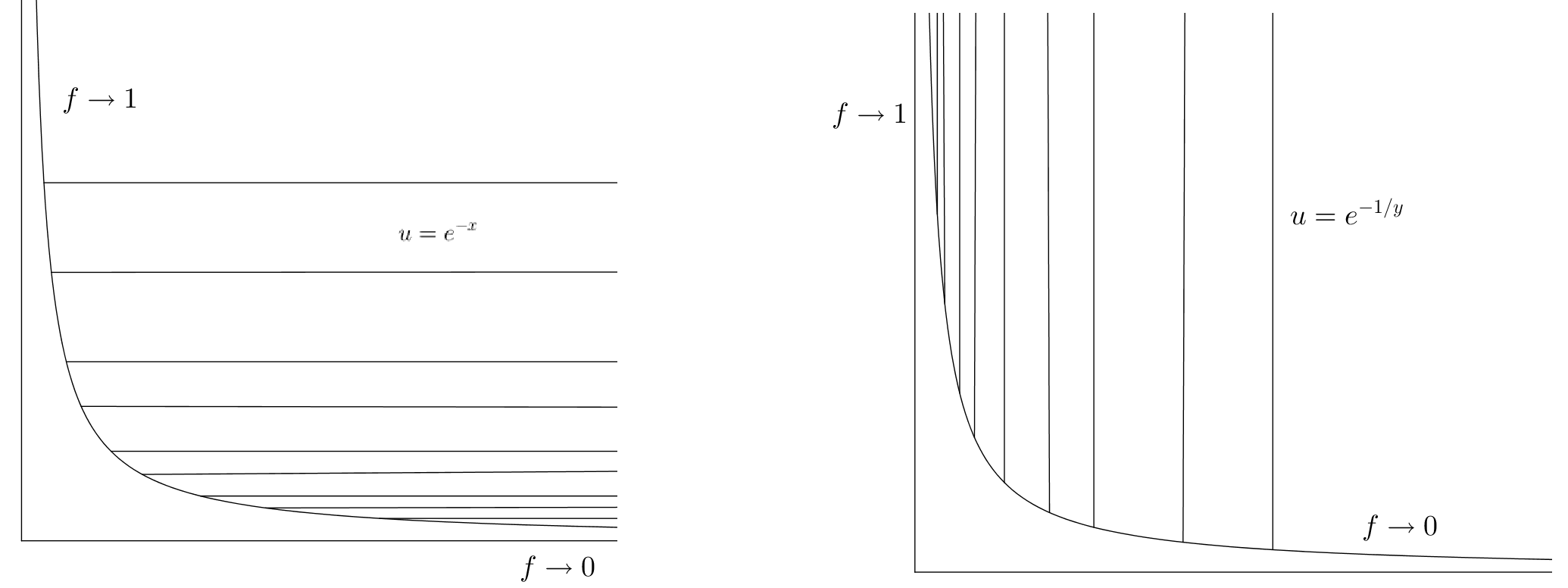}
    \caption{Nonuniqueness of minimisers for data in $C_b(\partial\Omega)$}
    \label{fig:nonuniqueness}
\end{figure}
\end{prz}

The reason why uniqueness fails in the example above is that the set $\Omega$ is not one-dimensional, but contains a cone. Therefore the level lines could ''choose'' a direction in which the solution propagates; this is the primary motivation for the considerations in the next subsection. Due to this phenomenon, the above example implies that in the unbounded case the structure theorems for least gradient functions, see \cite[Theorem 1.1]{Mor} and \cite[Theorem 1.1]{Gor1}, cannot hold - there exist multiple continuous functions which solve the least gradient problem and their derivatives do not agree on a set of positive Lebesgue measure.

\begin{uw}
As we rely in the uniqueness proof on the Sternberg-Williams-Ziemer construction, it generalises to the anisotropic setting in the same cases as in the continuity proof - when $\Omega \subset \mathbb{R}^2$ and $\phi$ is a norm such that $B_\phi(0,1)$ is strictly convex and for weighted least gradient problem, when $\phi(x, Du) = a(x)|Du|$ with $a \in C^2(\overline{\Omega})$. The reason is that uniqueness provided by the Sternberg-Williams-Ziemer construction stems from the maximum principle for minimal surfaces.
\end{uw}

\subsection{The case when $\Omega = \mathbb{R}^N$} 

When we consider unbounded domains, in previous subsections we only discussed boundary conditions on $\partial\Omega$, disregarding the limit behaviour. In this subsection, we are interested in the case when we impose the limit behaviour at infinity. We consider the following problem: given a function $f \in L^1(\partial B(0,1))$, we want to find a least gradient function $u \in BV_{loc}(\mathbb{R}^N)$ on $\mathbb{R}^N$ such that for $\mathcal{H}^{N-1}$-a.e. $x \in \partial B(0,1)$
$$ u(tx) \rightarrow f(x) \qquad \text{ as } t \rightarrow \infty.$$
Let us first focus on $\mathbb{R}^2$.

\begin{stw}
Let $\Omega = \mathbb{R}^2$. Denote by $\Gamma_1, \Gamma_2$ any two non-intersecting half-circles in $\partial B(0,1)$. Then all the least gradient functions in $\mathbb{R}^2$ are one-dimensional. In particular, the only limit values of $u$ at infinity are of the form $f = a \chi_{\Gamma_1} + b \chi_{\Gamma_2}$.
\end{stw}

\begin{dd}
Suppose that $u \in BV_{loc}(\mathbb{R}^2)$ is a function of least gradient. Fix any $t \in \mathbb{R}$. By Theorem \ref{tw:bgg} the function $\chi_{\{ u \geq t \}}$ is also a function of least gradient. Using the regularity theory for minimal sets proved by Giusti in \cite{Giu}, we obtain that each connected component of $\partial \{ u \geq t \}$ is in fact a smooth minimal surface. In two dimensions, it means that $\partial \{ u \geq t \}$ is a union of at most countably many parallel lines. However, we easily see that this union contains only one element - if it contained more than one element, we could replace any two of them with two U-shaped polygonal chains, which would have locally smaller total variations (analogously to what is presented on Figure \ref{fig:escape}).

Hence for each $t \in \mathbb{R}$ we have that $\{ u \geq t \}$ is either empty of $\{ u \geq t \} = l_t$. Furthermore, as for $t > s$ we have $\{ u \geq t \} \subset \{ u \geq s \}$, all of these lines are parallel to a line $l$ passing through the origin. Then $u$ is a function of one variable $z$, defined along the line $l'$ passing through the origin and perpendicular to $l$. The orientation of the line $l'$ is chosen so that $z$ is growing in the direction of angles from the interval $[\frac{\pi}{2}, \frac{3\pi}{2})$.


Let us denote the angle of incidence of $l$ by $\alpha_0$. Let $u(z) \rightarrow a$ as $z \rightarrow \infty$ and $u(z) \rightarrow b$ as $z \rightarrow - \infty$. Then, as for $\alpha \neq \alpha_0$ and $\alpha \neq \alpha_0 + \pi$ any ray starting from the origin has to intersect all the lines $l_t$, the limit value of $u$ at infinity is
$$ f(x,y) = f(\alpha) = \twopartdef{a}{\text{if } \alpha \in (\alpha_0, \alpha_0 + \pi)}{b}{\text{if } \alpha \in (\alpha_0 + \pi, \alpha_0 + 2\pi).}$$
Note that here, with a slight abuse of notation, at most one of $a, b$ may be equal to $+ \infty$ and at most one of $a, b$ may be equal to $- \infty$. \qed
\end{dd}

For a discussion in $\mathbb{R}^3$ we will need an additional result. It comes from the classical theory of minimal surfaces and is called the strong halfspace theorem, proved by Hoffman and Meeks in \cite{HM}. 

\begin{tw}\label{halfspace}
(\cite[Theorem 2]{HM}) Two proper, possibly branched, connected minimal surfaces in $\mathbb{R}^3$ must intersect, unless they are parallel planes.
\end{tw}

\begin{stw}
Let $\Omega = \mathbb{R}^3$. Denote by $\Gamma_1, \Gamma_2$ any two non-intersecting half-spheres in $\partial B(0,1)$. Then all the least gradient functions in $\mathbb{R}^3$ are either one-dimensional or have a single jump along a smooth minimal surface $S$. In particular, the only limit values of $u$ at infinity are of the form $f = a \chi_{\Gamma_1} + b \chi_{\Gamma_2}$.
\end{stw}

\begin{dd}
Suppose that $u \in BV_{loc}(\mathbb{R}^3)$ is a function of least gradient. Fix any $t \in \mathbb{R}$. As in the proof of the previous Proposition, we obtain that each connected component of $\partial \{ u \geq t \}$ is a smooth properly embedded minimal surface. By Theorem \ref{halfspace} there is only one connected component of $\{ u \geq t \}$ unless $\{ u \geq t \}$ is a union parallel planes; however, we argue as in the proof of the previous Proposition that in that case $\{ u \geq t \}$ is a single plane.

Now, consider $t > s$. As $\{ u \geq t \} \subset \{ u \geq s \}$, if $\partial \{ u \geq t \} \cap \partial \{ u \geq s\}$, then by Proposition \ref{swz:zasadamaksimum} these sets agree. Therefore, if for any $t \in \mathbb{R}$ the surface $S = \partial \{ u \geq t \}$ is not a plane, then for any $s \in \mathbb{R}$, by Theorem \ref{halfspace} $\partial \{ u \geq s \}$ intersects $S$ and therefore agrees with it; hence the function $u$ has a single jump across $S$ and is constant in $\mathbb{R}^3 \backslash S$. 

Now, we discuss the limit behaviour at infinity. If all sets of the form $\partial \{ u \geq t \}$ are planes, then they are parallel to a plane $\Pi$ passing through the origin and we proceed as in the proof of the previous Proposition to get that $f$ is a function with two values which are obtained on two disjoint halfspheres. If the function $u$ has only two values and a single jump across a smooth properly embedded minimal surface $S$, then recall that $S$ admits a limit tangent plane at infinity (for instance in the sense of \cite[Definition 6.1]{MP}) denoted by $\Pi$ (by definition, it passes through the origin); let us take any ray starting from the origin such that its direction does not lie in $\Pi$. Then the value of $u$ along that ray stabilises. Thus, the limit value of $u$ at infinity is
$$ f(x,y,z) = \twopartdef{a}{\text{in a halfsphere above } \Pi}{b}{\text{in a halfsphere below } \Pi.}$$
Again, with a slight abuse of notation at most one of $a, b$ may be equal to $+ \infty$ and at most one of $a, b$ may be equal to $- \infty$. \qed
\end{dd}

In higher dimensions the situation is much less clear. There are two main reasons for this: firstly, we do not have the halfspace theorem and therefore there may exist continuous least gradient functions which are not one-dimensional. Moreover, in dimensions eight and above there exist least gradient functions such that $\partial \{u \geq t \}$ may have singularities. To illustrate this, we recall the construction developed in \cite{BGG}.

\begin{prz}
(1) Let $\Omega = \mathbb{R}^8$. Let $C \subset \mathbb{R}^8$ denote the interior of the Simons cone, namely
$$ C = \{ x_1^2 + x_2^2 + x_3^2 + x_4^2 \geq x_5^2 + x_6^2 + x_7^2 + x_8^2 \}.$$
Then $u = \chi_C$ is a function of least gradient. However, the limit values of $u$ at infinity equal $f = \chi_{C \cap \partial B(0,1)}$, which is not constant on halfspheres. Moreover, the authors of \cite{BGG} construct (in the proof of Theorem A) a continuous function $F$ of least gradient such that the Simons cone is its zero level set. In particular, this function is not one-dimensional neither has any jumps.

(2) Let $N \geq 9$. Then the answer to the Bernstein problem is positive and there exist entire complete analytic minimal graphs in $\mathbb{R}^N$ which are not hyperplanes. Hence, we can construct a function of least gradient such that its level sets are translations of a single Bernstein graph; in particular, this function is not one-dimensional and we can ensure that is has no jumps.
\end{prz}

Overall, these results suggest that the formulation of the Dirichlet problem with boundary conditions at infinity is not the proper way to introduce the least gradient problem on unbounded domains, due to the fact that the problem in this formulation does not have any solutions in low dimensions, unless the prescribed data has a very specific form. Hence, the more proper approach is to restrict ourselves to strictly convex sets $\Omega$ which are not equal to the whole of $\mathbb{R}^N$ and consider the Dirichlet boundary data only on $\partial\Omega$. Furthermore, there is no need to impose boundary conditions at infinity for 
such sets, as the limit behaviour is regulated by the Dirichlet data on $\partial\Omega$.

{\bf Acknowledgements.} I would like to thank my PhD advisor, Piotr Rybka, for many fruitful discussions about this paper. This work was partly supported by the research project no. 2017/27/N/ST1/02418, "Anisotropic least gradient problem", funded by the National Science Centre, Poland.

\bibliographystyle{plain}
\bibliography{WG-references}

\end{document}